\documentclass[12pt,draftcls,onecolumn]{IEEEtran}
\usepackage{graphicx}      % include this line if your document contains figures
\usepackage{mathtools, verbatim}
\usepackage{amsmath, tabularx}
\usepackage{amsfonts}
\usepackage{amssymb}
\usepackage{amsbsy}
\usepackage{algorithm}
\usepackage{algorithmic}
\usepackage{enumerate}

\usepackage{enumitem}
\usepackage{color}

\usepackage{latexsym}
\usepackage{cite}
\newcommand{\eg}{{\it e.g.}}
\newcommand{\ie}{{\it i.e.}}

\newcommand{\reals}{\mathbb{R}}
\newcommand{\natur}{\mathbb{N}}

\newcommand{\zinf}{\boldsymbol{z}}
\newcommand{\xinf}{\boldsymbol{x}}
\newcommand{\uinf}{\boldsymbol{u}}
\newcommand{\yinf}{\boldsymbol{y}}

\newcommand{\einf}{\boldsymbol{e}}

\newcommand{\lambdainf}{\boldsymbol\lambda}

\newcommand{\nuinf}{\boldsymbol\nu}

\newcommand{\Hinf}{\boldsymbol H}

\newcommand{\Cinf}{\boldsymbol C}
\newcommand{\cinf}{\boldsymbol c}
\newcommand{\hinf}{\boldsymbol h}

\newcommand{\Ainf}{\boldsymbol A}
\newcommand{\Binf}{\boldsymbol B}
\newcommand{\Ginf}{\boldsymbol G}
\newcommand{\cxinf}{\boldsymbol {c_x}}

\newcommand{\Cxinf}{\boldsymbol {C_x}}

\newcommand{\Rinf}{\boldsymbol R}
\newcommand{\Qinf}{\boldsymbol Q}
\newcommand{\Winf}{\boldsymbol W}

\newcommand{\diag}{\mathop{\bf diag}}

\newcommand{\argmin}{\mathop{\rm argmin}}

\newcommand{\interior}{\mathop{\bf int}}

\newcommand{\prox}{\mathop{\bf prox}}

\newcommand{\tr}{{\top}}

\newcommand{\Hz}{\mathcal{H}_{\boldsymbol z}}
\newcommand{\HI}{\mathcal{H}_{\mathbf 1}}
\newcommand{\HII}{\mathcal{H}_{\mathbf 2}}
\newcommand{\Hu}{\mathcal{H}_{\boldsymbol u}}
\newcommand{\Hx}{\mathcal{H}_{\boldsymbol x}}
\newcommand{\Hl}{\mathcal{H}_{\boldsymbol \lambda}}

\newcommand{\lw}{l^2_{w}}

\newcommand{\comments}[1]{} 
\newcommand{\LeftEqNo}{\let\veqno\@@leqno}

\newtheorem{mydef}{Definition}
\newtheorem{mylem}{Lemma}
\newtheorem{myas}{Assumption}
\newtheorem{myrem}{Remark}

\newtheorem{thm}{Theorem}

\begin{document}

%\title{Solving the infinite-horizon constrained LQR problem using splitting techniques}

%\thanks[footnoteinfo]{Sponsor and financial support acknowledgment
%goes here. Paper titles should be written in uppercase and lowercase
%letters, not all uppercase.}

%\author[First]{Giorgos Stathopoulos} 
%\author[First]{Milan Korda} 
\title{Solving the Infinite-horizon Constrained LQR Problem using Accelerated Dual Proximal Methods}
\author{Giorgos Stathopoulos, Milan Korda and Colin N. Jones \thanks{G. Stathopoulos, M. Korda and C.N. Jones are with the Laboratoire d'Automatique, \'Ecole Polytechnique F\'ed\'erale de Lausanne (EPFL), Lausanne, Switzerland. E-mail: \{georgios.stathopoulos,milan.korda,colin.jones\}@epfl.ch .}}

\author{Giorgos~Stathopoulos,
        Milan~Korda,
        and~Colin~N.~Jones.% <-this % stops a space
\thanks{The authors are with the the Laboratoire d'Automatique, \'Ecole Polytechnique F\'ed\'erale de Lausanne (EPFL), Lausanne, Switzerland. E-mail: \{georgios.stathopoulos,milan.korda,colin.jones\}@epfl.ch .}}
%\thanks{Manuscript received April 19, 2005; revised December 27, 2012.}}

%\markboth{IEEE Transactions On Automatic Control, Vol. XX, No. Y, Month 2014}
%{Murray and Balemi: Using the style file IEEEtran.sty}

\maketitle
%\thispagestyle{plain}\pagestyle{plain}

%\begingroup
%\centering
%{\LARGE \textbf{Solving the infinite-horizon constrained LQR problem using splitting techniques} \\[1.5em]
%\large Giorgos Stathopoulos\Mark{1}, Milan Korda\Mark{1} and Colin N. Jones\Mark{1}}\\[1em]
%\begin{tabular}{*{1}{>{\centering}p{.60\textwidth}}}
%\Mark{1}Laboratoire d'Automatique\tabularnewline
%\'Ecole Polytechnique F\'ed\'erale de Lausanne\tabularnewline
%Station 9, CH-1015, Lausanne, Switzerland\tabularnewline
%\{georgios.stathopoulos,milan.korda,colin.jones\}@epfl.ch
%\end{tabular}\par
%\endgroup
%\maketitle

%\address[First]{Laboratoire d'Automatique, \'Ecole Polytechnique F\'ed\'erale de Lausanne, Station 9,
%CH-1015, Lausanne, Switzerland;\\ email:  \{georgios.stathopoulos,milan.korda,colin.jones\}@epfl.ch}

\begin{abstract}      
This work presents an algorithmic scheme for solving the infinite-time constrained linear quadratic regulation problem. We employ an accelerated version of a popular proximal gradient scheme, commonly known as the Forward-Backward Splitting (FBS), and prove its convergence to the optimal solution in our infinite-dimensional setting. Each iteration of the algorithm requires only finite memory, is computationally cheap, and makes no use of terminal invariant sets; hence, the algorithm can be applied to systems of very large dimensions. The acceleration brings in  `optimal' convergence rates $O(1/k^2)$ for function values and $O(1/k)$ for primal iterates and renders the proposed method a practical alternative to model predictive control schemes for setpoint tracking. In addition, for the case when the true system is subject to disturbances or modelling errors, we propose an efficient warm-starting procedure, which significantly reduces the number of iterations when the algorithm is applied in closed-loop. Numerical examples demonstrate the approach.

%The regulation problem is cast in its condensed formulation, \ie, written in terms of the input variables, %and subsequently recast in its dual domain. In this way, the original polytopic constraints on states and %inputs are transformed into nonpositivity constraints on the dual variables. The resulting infinite %dimensional problem admits a finite-dimensional representation provided the dual sequenced is suitably %initialized. 
    
% An algorithmic scheme for solving the constrained linear quadratic regulation problem in real time is presented in this work. We recruit a popular proximal gradient scheme, commonly known as the Forward-Backward Splitting (FBS). The regulation problem is cast in its condensed formulation, \ie, written in terms of the input variables, and subsequently recast in its dual domain.  In this way, the original polytopic constraints on states and inputs are transformed into nonpositivity constraints on the dual variables. The resulting infinite dimensional problem can be tackled in finite dimensions by observing that the dual sequence always vanishes after some (finite) time epoch, a fact that allows for a finite representation of the state and input sequences. The proposed approach makes no use of terminal invariant sets and generates the infinite horizon optimal solution to the problem, along with a sufficient horizon length for achieving it. In addition, a provably convergent accelerated version of the (typically slow) FBS is employed, rendering the proposed method a practical alternative to model predictive control for regulation.
\end{abstract}


\IEEEpeerreviewmaketitle
%===============================================================================

\section{Introduction}
An important extension of the famous result of \cite{Kalman} on the closed form solution of the infinite-horizon linear quadratic regulation (LQR) problem is the case where input and state variables are constrained. This problem is computationally significantly more difficult and has been by and large addressed only approximately. A prime example of an approximation scheme is \emph{model predictive control (MPC)}, which approximates the infinite-time constrained problem by a finite-time one. Stability of such MPC controllers is then typically enforced by adding a suitable terminal constraint and a terminal penalty. The inclusion of a terminal constraint limits the feasible region of the MPC controller, and, consequently, the region of attraction of the closed-loop system. In practical applications, this problem is typically overcome by simply choosing a ``sufficiently'' long horizon based on process insight (e.g., dominant time constant). Closed-loop behavior is then analyzed a posteriori, for instance by exhaustive simulation.

There have been few results addressing directly the \emph{infinite-horizon constrained LQR (CLQR)} problem. Among the most well-known efforts are the works \cite{Chmielewski1996121}, \cite{CLQR} and \cite{Grieder2004701}. The authors of \cite{Chmielewski1996121} suggest a scheme for offline computation of a sufficient horizon length. The solution of the corresponding quadratic program (QP) is then equivalent to the original infinite dimensional problem. The reported results are somewhat conservative, while the offline part of the proposed algorithm can be computationally prohibitive since it involves the solution of a possibly nonconvex problem, the computation of a positively invariant set as well as a vertex enumeration problem. The authors of \cite{CLQR} extend the work of \cite{Sznaier} by solving a sequence of QPs of finite horizon length, which is monotonically non-decreasing. After each QP has been solved, the inclusion of the final state in a positively invariant set associated to the optimal unconstrained LQ controller is checked; if the final state is not included in the set, the horizon was insufficient and has to be increased. Finally, the authors in \cite{Grieder2004701} employ parametric quadratic programming and a reachability analysis approach to compute the least conservative horizon length that ensures optimal infinite horizon performance. Although this work provides the exact necessary horizon length for the feasible set of initial states, it suffers from tractability issues since it is based on state-space partitioning and thus can only be applied to small systems.

Our approach is inspired by the framework of \emph{proximal gradient methods}, a class of algorithms known for their simplicity when it comes to solving convex optimization problems with simple constraints~\cite{parikh2014proximal}. From this family of algorithms, we use the \emph{(accelerated) version of the Forward-Backward Splitting method ((A)FBS)} (see, \eg~\cite[Chapter~25]{book_comb} and \cite{AFBS_pock}). The idea is to condense the problem and describe it in terms of its input variables, then formulate the dual of the corresponding (infinite-dimensional) QP, and apply the (A)FBS method to solve it. More specifically, the method decomposes the QP into two subproblems, the first one being an infinite-dimensional least squares problem and the second one a simple clipping of an infinite sequence to the non-positive orthant. The subproblems are solved repeatedly (with the solution of one influencing the cost function of the other) until convergence to the solution of the original problem. This is in contrast to the approach of~\cite{CLQR}, which requires the solution of a sequence of \emph{constrained} QPs. We show that both sub-problems of the proposed algorithm can be solved tractably (which is not a priori obvious since we are working with infinite sequences), the first one by solving a finite-dimensional system of linear equations (with the possibility to pre-factorize the matrices) and the second one by simple clipping of finitely many real numbers on the non-positive real line. Convergence of the scheme (with rate $O(1/k^2)$ for function values and $O(1/k)$ for primal iterates) to the optimal infinite-horizon sequence is guaranteed under mild assumptions. Therefore the proposed algorithmic scheme provides a means to compute the solution of the infinite-horizon constrained LQR problem with guaranteed convergence.

This work is based on our recent result in~\cite{Sta_Kor_Jon_CLQR}, where the same problem is solved by employing a primal-dual scheme, the \emph{Alternating Minimization Algorithm (AMA)}~\cite{Tseng:1991}, which is equivalent to FBS. The approach is equivalent to splitting the infinite-horizon constrained LQR problem into an \emph{unconstrained} LQR problem and a proximal minimization problem. In the current work, the result is developed in both theoretical and practical directions. From the theoretical perspective, we (i) drop the limiting assumption for positive-definiteness of the state penalty matrix $Q$ that existed in~\cite{Sta_Kor_Jon_CLQR}, and, more importantly, we (ii) prove convergence of an accelerated counterpart of the original FBS method, a fact that enables a (significantly) faster convergence rate. From the practical perspective, we provide a fully implementable method, competitive for real-time control. We (i) eliminate the need for knowledge of uncomputable quantities, (ii) propose computationally efficient ways to solve the optimization subproblems and (iii) propose a warm-starting scheme that performs well in practice.

The paper is organized as follows: In Section \ref{sec:intro} we introduce the problem and express it in terms of its dual variables by means of the proximal splitting framework. In Section~\ref{sec:afbs} we present the accelerated forward backward splitting algorithm to solve the problem and show that each iteration of the algorithm can be carried out tractably. The convergence proofs for this scheme are given in Section~\ref{sec:conv}, while Section~\ref{sec:lin_sys_solve} discusses the computational aspects that make the algorithm practical to use. Section~\ref{sec:examples} presents two numerical examples: A toy example of an unstable system with two states and one control input illustrates the main features of the algorithm. Subsequently, we demonstrate the practical applicability of the algorithm on a linearized model of a quadcopter with 12 states, 4 inputs and polytopic constraints. Finally, Appendices~\ref{App.A},~\ref{App.B} and~\ref{App.C} provide the proofs for the results presented in Section~\ref{sec:conv}.

%\section{Notation}{\label{sec:notation}}

\section{Problem statement and dualization}{\label{sec:intro}}
The constrained regulation problem can be written in the general form
\begin{equation}
\label{mpc-initial_reg}
\begin{array}{ll}
\mbox{minimize} & \frac{1}{2}\overset{\infty}{\underset{i=0}{\sum}} x_i^\tr Qx_i + u_i^\tr Ru_i\\
\mbox{subject to} & x_{i+1} = A x_i + B u_i , \quad i\in\mathbb{N}\\
& x_0 = x_\mathrm{init}\\
& C_xx_i \leq c_x\\
& C_uu_i \leq c_u\enspace,
\end{array}
\end{equation}
with variables $x_i\in \mathbb{R}^n$ and $u_i\in \mathbb{R}^m$, and data $c_x\in \mathbb{R}^{p_x}$ and $c_u\in \mathbb{R}^{p_u}$. 

We make the following standing assumptions:
\begin{myas}\label{as:basic_reg}
The pair $(A,B)$ is stabilizable, the optimal value of problem~(\ref{mpc-initial_reg}) is finite, the set \[\mathcal{X} := \{x\in\mathbb{R}^n\mid C_xx \le c_x\}\] contains the origin in its interior, the matrix $C_u$ has full column rank, the matrix $Q$ is positive semidefinite and $R$ is positive definite.
\end{myas}
\begin{myrem}[Stability]
Clearly, under Assumption~\ref{as:basic_reg}, if the problem~(\ref{mpc-initial_reg}) is feasible and the pair $(A,\sqrt{Q})$ detectable, then the control sequence optimal in~(\ref{mpc-initial_reg}) is stabilizing. Therefore, there is no need to enforce stability in an ad hoc way as is commonly done when the infinite-time problem~(\ref{mpc-initial_reg}) is approximated by a finite-time one solved in a receding horizon fashion.
\end{myrem}

From now on we write infinite sequences and infinite-dimensional operators in bold font.

The problem can be rewritten in the dense form, \ie, by writing the states as functions of the inputs. This is done by defining the operators
\[
\Ainf = \left[\begin{array}{c}
A \\ A^2\\ \vdots\end{array}\right],\quad
\Binf = \left[\begin{array}{ccccc}
B&0&0&\cdots \\ AB&B&0&\cdots\\A^2B&AB&B&\cdots \\ \vdots&\vdots&\vdots&\ddots \end{array}\right]\enspace,
\]
\[
\Qinf = \diag(Q,Q,\ldots)\enspace, \quad \Rinf = \diag(R,R,\ldots)\enspace,
\]
\[
\Cxinf = \diag(C_x,C_x,\ldots)\enspace, \quad \cxinf = (c_x,c_x,\ldots)\enspace,
\]
\[
\Cinf_i = \left[\begin{array}{c}\einf_i\otimes C_u \\ \left[\Cxinf\Binf\right]_i\end{array}\right], \quad c_i = \left[\begin{array}{c} c_u \\ c_x -C_x A^i x_{\mathrm{init}}\end{array}\right]\enspace,
\]
\[
\Hinf = \Binf^{\star}\Qinf\Binf+\Rinf, \quad \Ginf = \Ainf^{\star}\Qinf\Binf, \quad \bar{F} = \Ainf^{\star}\Qinf\Ainf+Q\enspace,
\]
\begin{equation}{\label{operators}}
\Cinf = \left[\begin{array}{ccc}\Cinf_1^{\star} & \Cinf_2^{\star} & \ldots \end{array}\right]^{\star}, \quad \cinf = [c_1^\star, c_2^\star,\ldots]^\star \enspace,
\end{equation}
where we denote by $\einf_i$ the (infinite dimensional) row vector with only one nonzero element equal to one at position $i$, by $\left[\cdot\right]_i$ the $i^{\mathrm{th}}$ block row of size $p_x\times \infty$ of the corresponding operator and by $^\star$ the adjoint of an operator (i.e., the infinite-dimensional analogue of transpose; see Appendix~A for a brief introduction to operators). 
Using the above, (\ref{mpc-initial_reg}) can be written in the form
\begin{equation}
\label{eq:generic_pr}
\begin{array}{ll}
\mbox{minimize} & \frac{1}{2}  \uinf^\star \Hinf\uinf  +  \hinf^\star \uinf  + r\\
\mbox{subject to} & \Cinf\uinf \leq \cinf\enspace,
\end{array}
\end{equation}
where $\hinf=\Ginf^{\star}x_{\mathrm{init}}$, $r=\frac{1}{2}  x_{\mathrm{init}}^\star \bar{F} x_{\mathrm{init}} $, $\uinf$ is the infinite sequence
$\uinf := [u_0^\star,u_1^\star,\ldots]^\star$ and $\Hinf\colon\Hu\to\Hu$, $\Cinf\colon\Hl\to\Hu$, $\cinf\in\Hl$, where $\mathcal{H}$'s are suitable Hilbert spaces specified next: Any sequence
\[
\zinf := (z_0,z_1,\ldots)
\]
is viewed as an element of an $l^2$-weighted (or $\lw$) real Hilbert space~$\Hz$ (see Appendix~A, Definition~\ref{def:wl2_hilbert_space}) induced by the inner product
\begin{equation}{\label{def::l2_space}}
\langle \zinf,\yinf\rangle  = \sum_{i=0}^\infty w^iz_i^\tr y_i\;,\quad \forall\, \yinf \in \Hz,\, \zinf \in \Hz\enspace,
\end{equation}
where $w$ is an appropriately chosen weight (see Appendix~A, Definition~\ref{def:wl2_hilbert_space}). The norm of any $\zinf \in \Hz$ is thus given by
\[
  \| \zinf  \|_{\Hz} := \sqrt{\langle \zinf,\zinf\rangle }   = \sqrt{\sum_{i=0}^\infty w^i \|z_i\|_2^2}\;\enspace.
\]
Unless stated otherwise, for the rest of the paper by a Hilbert space we mean the $\lw$ real Hilbert space as just introduced. 

The dual of problem (\ref{eq:generic_pr}) can be derived by making use of Legendre-Fenchel duality (see, \eg, Chapter 7, \cite{bertsekas03}).
By defining 
\begin{equation}{\label{eq:f_g_def}}
f(\uinf) = \frac{1}{2}  \uinf^\star \Hinf\uinf  +  \hinf^\star \uinf  + r, \quad g(\Cinf\uinf-\cinf) = \delta_-(\Cinf\uinf-\cinf)\enspace
\end{equation}
and $\delta_-(\cdot)$ being the indicator function for the nonpositive orthant, we can rewrite (\ref{eq:generic_pr}) as
\[
\underset{\uinf}{\min}\left\{f(\uinf)+g(\Cinf \uinf-\cinf)\right\}
\]
and then express the Lagrange dual problem by using the dual variable $\lambdainf$ as
\begin{align*}
& \underset{\lambdainf}{\max}\left\{\underset{\uinf}{\min}\left\{f(\uinf)+\langle \cinf-\Cinf \uinf,\lambdainf\rangle \right\}-g(\lambdainf)\right\} \Leftrightarrow \\
& \underset{\lambdainf}{\max}\Big\{-\underset{\uinf}{\max}\left\{\langle \Cinf^{\star}\lambdainf,\uinf\rangle -f(\uinf)\right\}-g(\lambdainf)+\langle \lambdainf,\cinf\rangle \Big\} \Leftrightarrow \\
& \underset{\lambdainf}{\max}\left\{-f^{\star}(\Cinf^{\star}\lambdainf) - g(\lambdainf) + \langle \lambdainf,\cinf\rangle \right\}\enspace,
\end{align*}

The solution can be derived by solving
\begin{equation}{\label{app:dual}}
 \underset{\lambdainf}{\min}\left\{f^{\star}(\Cinf^{\star}\lambdainf) + g(\lambdainf) - \langle \lambdainf,\cinf\rangle \right\}\enspace.
\end{equation}
The function involved in the minimization problem (\ref{app:dual}) is
\begin{equation}\label{eq:fstar}
f^{\star}(\Cinf^{\star}\lambdainf) = \frac{1}{2}\langle \Cinf^{\star}\lambdainf,\Hinf^{-1}\Cinf^{\star}\lambdainf\rangle 
                                   - \langle \Cinf^{\star}\lambdainf,\Hinf^{-1}\hinf\rangle + \frac{1}{2}\langle \hinf,\Hinf^{-1}\hinf\rangle  - r\enspace.
\end{equation}
\noindent
Finally, the problem we are interested in solving can be cast in a more compact form as:
\begin{equation}
\label{eq:generic_dual}
\begin{array}{ll}
\mbox{minimize} & F(\lambdainf) := h^{\star}(\lambdainf)+ \delta_-(\lambdainf),
\end{array}
\end{equation}
with variables $\lambdainf$, $h^{\star}(\lambdainf) = f^{\star}(\Cinf^{\star}\lambdainf)- \langle \lambdainf,\cinf\rangle  $ and $h^{\star}(\lambdainf),\delta_-(\lambdainf)$ being proper lower semi-continuous convex functions from $\Hl$ to $\reals\cup \left\{+\infty\right\}$.

Before proceeding, we elaborate on the reasons why the original problem~(\ref{mpc-initial_reg}) had to be reformulated in order to solve it. There are two reformulations, namely, posing the problem as a function of the input sequences only, resulting in (\ref{eq:generic_pr}), and dualization of (\ref{eq:generic_pr}), resulting in (\ref{eq:generic_dual}). The reason for considering the condensed formulation is the need for \emph{strong convexity} of the primal objective, which implies the Lipschitz continuity of $\nabla h^{\star}(\cdot)$~\cite[Corollary~18.16]{book_comb}. By using the condensed form we avoid the restrictive assumption of $Q\succ 0$ required in~\cite{Sta_Kor_Jon_CLQR}.

The reason for considering the dual problem is simplicity in the evaluation of the \emph{proximal operator} (see Appendix~\ref{App.A}, Defintion~\ref{def:proximal}) of the function $\delta_-(\cdot)$, which is a simple projection on the non-positive orthant (i.e., componentwise clipping) as opposed to the primal case where one would have to project on a generic polytope of the form $\Cinf\uinf \leq \cinf$.

Both steps are crucial for successfully applying the \emph{Forward-Backward Splitting (FBS)} method, as presented and analyzed in~\cite[Theorem~25.8]{book_comb}

\section{Solution using AFBS}\label{sec:afbs}

Problem~(\ref{eq:generic_dual}) is a composite minimization problem (i.e., a minimization of a sum of a smooth and a non-smooth function) and will be a solved using an \emph{accelerated forward-backward splitting (AFBS)} method, which is an accelerated variant of the forward-backward splitting proposed in~\cite[Theorem~25.8]{book_comb}. The acceleration comes from a Nesterov-like momentum sequence from~\cite{AFBS_chambolle}. 

This modification allows acceleration of the FBS method in a very simple manner and adds practically zero computational complexity per iteration, in comparison to the original version. The idea can be traced back to Polyak and the so called heavy ball methods~\cite{Polyak19641} for minimizing a smooth convex function 
$f(\cdot)$:
\begin{align*}%{\label{heavy_ball}}
 \hat{\lambda}^k     &= \lambda^k + \alpha^k(\lambda^k-\lambda^{k-1})\\
 \lambda^{k+1} &= \hat{\lambda}^k - \rho^k\nabla f(\lambda^k)\enspace, 
\end{align*}
where $\alpha^k\in[0,1)$ is an extrapolation factor and $\rho^k$ a stepsize parameter. This seemingly small change of updating the new iterate as a linear combination of the two previous iterates greatly improves the performance of the original gradient scheme. 

In his seminal paper~\cite{nest83}, Nesterov modified the heavy ball method by simply evaluating the gradient at the extrapolated point $\hat{\lambda}^k$ instead of $\lambda^k$. In addition, he proposed a special formula for computing the relaxation sequence $\left(\alpha^k\right)_{k\in\natur}$, resulting in an optimal convergence rate for the scheme. Subsequently, G\"{u}ler extended Nesterov's results to the proximal point algorithm, handling the minimization of the sum of a smooth and a nonsmooth function~\cite{Guler:1991:CPP:104219.104233}.

%The work~\cite{AFBS_chambolle} proposes a whole family of momentum sequences, including the Nesterov's original one, and analyzes them in a Hilbert-space setting as required here.

The acceleration that optimal first order methods enjoy comes solely from the additional momentum of the minimizer sequence due to the relaxation sequence $\alpha^k$. 
Among the many works that derived optimal relaxation sequences similar to that of Nesterov, a distinguished one is that of Beck and Teboulle in~\cite{fista}, who used Nesterov's framework to accelerate the FBS method. The resulting algorithm, commonly known as \emph{FISTA (Fast Iterative Shrinkage Thresholding Algorithm)} gives an optimal, $1/k^2$ convergence rate in terms of the function values.
%\[
%f(\lambda^k)-f(\lambda^\star) \le \frac{2L(f)}{(k+1)^2}\|\lambda^0-\lambda^\star\|^2\enspace.
%\]
%
This result has valuable practical and theoretical implications, however, convergence of the iterates of the FISTA algorithm has not yet been proven. Building upon FISTA, Lorenz and Pock suggest in \cite{AFBS_pock} a modified version of the method that achieves weak convergence of the iterates to the solution of the optimization problem in a Hilbert space setting, yet not proving anything about convergence of the function values.

The several `optimal' relaxation sequences have been put under a common framework in the recent work~\cite{AFBS_chambolle}. The authors showed that any sequence $\left(t^k\right)_{k\in\natur}$ of the form $t^k=\frac{k+a-1}{a}$, with $a > 2$ satisfies the inequality $(t^k)^2-t^k\le (t^{k-1})^2, \; k\ge 2$. Then the sequence defined as $\alpha^k=\frac{t^k-1}{t^{k+1}}$ allows for the optimal $O(1/k^2)$ convergence rate in terms of the function values, as well as weak convergence of the iterates, and hence provides an optimal first order scheme. We denote the algorithm emanating from this scheme the Accelerated Forward-Backward Splitting (AFBS) algorithm, and write it down for problem~(\ref{eq:generic_dual}) as follows:

%\begin{scriptsize}
\begin{algorithm}[H]
\caption{AFBS for Problem (\ref{eq:generic_dual})}
\label{al:AFBS}
\begin{algorithmic}  
\STATE 0: Initialize \scalebox{0.92}{$\lambdainf^0 = \mathbf{0}$, $a > 2$, $\alpha^0=0$}, \\
				\quad	 \scalebox{0.92}{$L(h^\star)-\mbox{ a Lipschitz constant of } \nabla h^\star$}.
\REPEAT
\STATE 1: $\alpha^k = \frac{k-1}{k+a}$,\; $k\ge 1$
\STATE 2: $\hat{\lambdainf}^{k} = \lambdainf^k + \alpha^k(\lambdainf^k-\lambdainf^{k-1})$
\STATE 3: $\lambdainf^{k+1} = \min\{\hat{\lambdainf}^k - \frac{1}{L(h^\star)}\nabla h^{\star}(\hat{\lambdainf}^k),0\}$
\UNTIL{termination condition is satisfied.}
\end{algorithmic}
\end{algorithm}
%\end{scriptsize}
The iterative scheme above is very simple: in this case the AFBS boils down to a variant of the fast projected gradient method as proposed by Nesterov~\cite{nesterov2004introductory}. In order to apply algorithm~(\ref{al:AFBS}), we need to to be able to
\begin{itemize}
\item evaluate the gradient of $h^\star(\cdot)$
\item represent $\lambdainf^k$ and $\hat{\lambdainf}^k$  using a finite amount of memory.
\end{itemize}
The remaining steps of the algorithm are simple scalar or vector updates  or componentwise clipping on the non-positive orthant (Step 3), both of which can be carried out inexpensively provided that $\lambdainf^k$ can be represented using finite amount of memory.

In the rest of the text $P_{\mathrm{LQ}}$ and $K_{\mathrm{LQ}}$ will denote the positive-semidefinite solution to the discrete-time algebraic Riccati equation and the corresponding LQ optimal state feedback matrix associated with the matrices $(A,B,Q,R)$.

We start by finding the gradient of $h^\star(\cdot)$. To this end, define the Lagrangian of the (infinite dimensional) problem~(\ref{eq:generic_pr}), written in its equivalent dense form and truncated at any $T\ge 0$, as
\begin{equation}\label{eq:L}
\mathcal{L}(\uinf, \xinf, \hat\lambdainf | T) \; = \; x_T^\top P_{\mathrm{LQ}}x_T + \sum_{i=0}^{T-1} x_i^\top Q x_i + u_i^\top R u_i - w^{i} \begin{bmatrix}C_x x_i - c_x\\ C_u u_i - c_u \end{bmatrix}^\top \hat\lambda_i. \nonumber
\end{equation}

\begin{mylem}\label{lem:grad}
If $\hat \lambdainf$ is such that $\hat \lambda_i = 0$ for all $i \ge T$, then
\begin{equation}\label{eq:grad}
[\nabla h^\star(\hat \lambdainf)]_i = 
w^{i}\begin{bmatrix}C_x \hat{x}_i - c_x\\ C_u \hat{u}_i - c_u \end{bmatrix},
\end{equation}
where
\begin{equation} \label{eq:Lmin}
(\hat{\xinf}, \hat{\uinf}) = \underset{\uinf,\xinf}{\argmin}\{ \mathcal{L}(\uinf, \xinf, \hat \lambdainf | T) \mid 
  x_{i+1} = Ax_i + Bu_i  ,\, x_0 = x_{\mathrm{init}}\}
\end{equation}
for $i \in \{0,\ldots, T \}$, and
\[
\hat{x}_i = (A+BK_{\mathrm{LQ}})^{i-T}\hat{x}_T,\quad \hat{u}_i = K_{\mathrm{LQ}}\hat{x}_i
\]
for $i > T$.
\end{mylem}
\begin{IEEEproof}
This is a standard result from duality (see, e.g., \cite[Chapter~5]{cvxbook}), noticing that the minimization in~(\ref{eq:Lmin}) is equivalent to the minimization of the Lagrangian of~(\ref{eq:generic_pr}) under the assumption that $\hat \lambda_i = 0$ for all $i \ge T$.
\end{IEEEproof}
Lemma~\ref{lem:grad} gives us a way to compute the gradient of $h^\star$. Clearly, this gradient is an infinite sequence and therefore cannot be stored directly, but it is available to us explicitly for $i \in \{0,\ldots T\}$ and implicitly, through the dynamics of the system $\hat x_+= (A+BK_{\mathrm{LQ}})\hat x$, for $i > T$.

Now we show that $\lambdainf^k$ and $\hat \lambdainf^k$ can be represented using a finite amount of memory.
\begin{mylem}
If $\lambdainf^k$ and $\hat \lambdainf^k$ are generated by Algorithm~\ref{al:AFBS}, then for each $k$ there exists a $T^k < \infty$ such that $\lambda_i^k = 0$ and $\hat{\lambda}_i^k = 0$ for all $i \ge T^k$.
\end{mylem}
\begin{IEEEproof}
We have $\lambda^0_i = 0$ and $\hat{\lambda}_i^0 = 0$ for all $i \ge 0$ and hence $T^0 = 0$. Assume now $k > 0$ and $\lambda_i^k = 0$ and $\hat{\lambda}_i^k = 0$ for all $ i\ge T^k $. Then, evaluating $\nabla h^\star(\lambdainf^k)$ using Lemma~\ref{lem:grad} with $T = T^k$, we see that the sequence $(\hat{x}_i,\hat{u}_i)$ is generated by the unconstrained LQ controller for $i > T^k$ and hence converges to the origin. Since the set $\mathcal{X}$ has the origin in the interior we conclude that there exists a time $\hat{T}^{k+1}<\infty$ such that $C_x\hat{x}_i \le c_x$ and $C_u\hat{u}_i \le c_u$ for all $i \ge \hat{T}^{k+1}$. We define $T^{k+1} = \max\{T^k, \hat{T}^{k+1}\}$. In view of~(\ref{eq:grad}), we conclude that $\hat{\lambdainf}^k - \frac{1}{L(h^\star)}\nabla h^{\star}(\hat{\lambdainf}^k) \ge 0$ for all $i \ge T^{k+1}$. Therefore, in view of Step 3 of Algorithm~\ref{al:AFBS}, we have $\lambda_i^{k+1} = 0$ for all $i \ge T^{k+1}$. Finally, $\hat{\lambdainf}^{k+1}$ is a linear combination of $\lambdainf^{k}$ and $\lambdainf^{k+1}$ and hence $\hat{\lambda}_i^{k+1}=0$ for all $i \ge T^{k+1}$. 
\end{IEEEproof}

To determine $T^{k+1}$ computationally (given $T^k$ and $\hat\xinf^{k+1}$ and $\hat\uinf^{k+1}$) we simply find the first time $T^\mathcal{S}$ that $\hat x^{k+1}_i$ enters a given subset $\mathcal{S}$, with $0\in\interior \mathcal{S}$, of the maximum positively invariant set of the system $\hat x_+ = (A+BK_{\mathrm{LQ}})\hat x$ subject to the constraint $\left[\begin{array}{c}C_uK\\C_x\end{array}\right]\hat x \leq \left[\begin{array}{c}c_u\\c_x\end{array}\right]$. The time $T^{k+1}$ is then equal to the first time no less than $T^k$ such that $C_x \hat x_{i+1}^{k+1}\le c_x$ and $C_u \hat u_i^{k+1} \le c_u$ simultaneously hold for all $i\in \{ T^{k+1}, \ldots, T^\mathcal{S} \}$. More formally, we have the equality
\begin{equation}\label{eq:Tk_update}
T^{k+1} = \min\big\{ T\ge T^{k} \mid \exists\, T^\mathcal{S} \; \mathrm{s.t.}\; C_x \hat x_{i+1}^{k+1}< c_x, 
                     C_u \hat u_i^{k+1} < c_u \;\forall\, i\in \{T,\ldots, T^\mathcal{S}\}, 
                     \mathrm{and} \;\; \hat x^{k+1}_{T^\mathcal{S}} \in \mathcal{S}   \big\}. 
\end{equation}
\begin{myrem}[Computation of $T^k$]\label{rem:sim}
In practice, to determine $T^{k+1}$ after solving~(\ref{eq:Lmin}), we iterate forward the system dynamics $\hat x_+ = (A+BK_{\mathrm{LQ}}) \hat x$ starting from the initial condition $\hat x^{k+1}_{T^k+1}$ until $\hat x^{k+1}_i \in \mathcal{S}$.
\end{myrem}
\begin{myrem}[Set $\mathcal S$]\label{rem:S}
The set $\mathcal{S}$ is determined offline and is \emph{not} required to be invariant. A good candidate is the set $\{x\mid x^\tr P_{\mathrm{LQ}} x \le 1\}$ scaled such that it is included in $\left\{x\mid \left[\begin{array}{c}C_uK\\C_x\end{array}\right]x \leq \left[\begin{array}{c}c_u\\c_x\end{array}\right]  \right\}$, or any subset of this set containing the origin in the interior.
\end{myrem}

Now we are ready to formulate an implementable version of the abstract algorithm~\ref{al:AFBS}:

\begin{scriptsize}
\begin{algorithm}[H]
\caption{AFBS for the CLQR problem}
\label{al:AFBS_CLQR}
\begin{algorithmic} 
\REQUIRE $ x_{\mathrm{init}}$, $Q\succeq 0,R\succ 0$, $C_u$ of full column rank, $w=\min\left\{1,\frac{1}{\lambda^2_{\mathrm{max}}(A)}\right\}$, $a > 2$, $\alpha^0=0$
\STATE a: Determine $P_{\mathrm{LQ}}$, $K_{\mathrm{LQ}}$ solving the unconstrained LQR
\STATE \hspace{4.2mm} problem associated with the matrices $(A,B,Q,R)$.
\STATE b: Determine a set $\mathcal{S}$, with $0\in\interior \mathcal{S}$, included in
\STATE \hspace{4.5mm} any positively invariant set for the system
\STATE \hspace{4.4mm} $x_+ = (A+BK_{\mathrm{LQ}})x$ subject to the constraint
\STATE \hspace{4.2mm} $\left[\begin{array}{c}C_uK\\C_x\end{array}\right]x \leq \left[\begin{array}{c}c_u\\c_x\end{array}\right]$. See Remark~\ref{rem:S}.
\STATE c: \scalebox{0.92}{\emph{Initialize $\lambdainf^0 = \mathbf{0}$, $T^0 = 0$, $a\ge 2$,}} \\
				\quad	 \scalebox{0.92}{\emph{$L^0>0 \mbox{ or } L(h^\star)-\mbox{ a Lipschitz constant of } \nabla h^\star$ (optional)}}.
\FOR{\texttt{$k=0,\ldots$}}
\STATE 1: $\alpha^k = \frac{k-1}{k+a}$,\; $k\ge 1$
\STATE 2: $\hat{\lambda}_i^{k} = \lambda^k_i + \alpha^k(\lambda^k_i-\lambda_i^{k-1}),\quad i=1,\ldots,T^k$
\STATE 3: Set \begin{align*}
				(\hat{\xinf}^{k+1}, \hat{\uinf}^{k+1}) &= \underset{\uinf,\xinf}{\argmin}\{
				\mathcal{L}(\uinf, \xinf, \hat \lambdainf^k | T^k) \mid 
				x_{i+1} = Ax_i + Bu_i  ,\, x_0 = x_{\mathrm{init}}\},\;
				i = 0,\ldots, T^k \\
				\hat{x}_{i+1}^{k+1} &= (A+BK_{\mathrm{LQ}}) \hat{x}_i^{k+1}, \quad i > T^k
			\end{align*}
\STATE 4: Determine $T^{k+1}$ (see Remark~\ref{rem:sim})
\STATE 6: Choose stepsize $\rho$ (see Remark~\ref{rem:stepize})
\STATE 7: Set $\lambda^{k+1}_i = \min\left(\lambda^k_i - \rho^{k+1}  w^{i}  \begin{bmatrix}C_x \hat{x}_i^{k+1} - c_x\\ C_u \hat{u}_i^{k+1} - c_u \end{bmatrix},0\right)$.
\STATE 8: If termination condition is satisfied, solve KKT system (see Remark~\ref{rem:termination})
\ENDFOR
\end{algorithmic}
\end{algorithm}\begin{footnotesize}
\end{footnotesize}
\end{scriptsize} 

\begin{myrem}[Stepsize]\label{rem:stepize}
In Step 6 of the algorithm, a stepsize is selected. One option is to fix a constant stepsize $\rho=\frac{1}{L(h^\star)}$, which needs a global Lipschitz constant of $\nabla h^\star (\cdot)$. This can be computed offline (see Section~\ref{sec:stepsize}). Alternatively, one can use a \emph{backtracking stepsize rule}, inspired from~\cite{fista}, which can be used in combination with the global estimate. The procedure is analyzed in Section~\ref{sec:lin_sys_solve}.
\end{myrem}

\begin{myrem}[Role of the weight $w$]
It is worth mentioning that working in the weighted Hilbert space $\lw$ is much more than a mathematical formalism and has serious practical implications. In the case of unstable systems, a nontrivial sequence of weights has to be chosen such that the growth of the largest unstable eigenvalue of the state matrix $A$ is bounded by a faster decaying sequence so that the operator $\Cinf$ remains bounded, an assumption necessary for applying the proposed method. At the same time the sequence of weights will act as a left and right preconditioner on the Hessian operator of the quadratic form in~(\ref{eq:generic_dual}), as well as a left preconditioner of $\Cinf$. This scaling can seriously affect the numerical performance of the proposed algorithms, as we will see in subsequent sections. Note that for stable systems the sequence can be trivially set to 1, hence no scaling occurs. These claims are explained in more detail in Appendix~\ref{App.C}.
\end{myrem}

\begin{myrem}[Termination]\label{rem:termination}
 Algorithm~\ref{al:AFBS_CLQR} terminates when a prespecified accuracy is reached in terms of the progress of the dual sequence. The extracted primal sequence is given by
 \begin{equation} \label{eq:primalSequence}
(\xinf^k, \uinf^k) = \underset{\uinf,\xinf}{\argmin}\{ \mathcal{L}(\uinf, \xinf, \lambdainf^k |  T^k) \mid
   x_{i+1} = Ax_i + Bu_i  ,\, x_0 = x_{\mathrm{init}}\}
\end{equation}
 for $i \in \{0,T^k\}$ and $x_i^k = (A+BK_{\mathrm{LQ}})^{i-T^k}x_{T^k}$ for $i > T^k$. In Theorem~\ref{theor:master} below it is proven that $(\xinf^k, \uinf^k)$ tends to the optimal constrained LQ solution. At any finite iterate, however, the sequence $(\xinf^k, \uinf^k)$ may violate the constraints. In order to remedy this we solve upon termination an equality-constrained QP where we minimize the objective function subject to the active constraints at optimality. The active constraints can be (approximately) detected by looking at the nonzero values of the dual vector $\lambdainf^k$ at termination. This step comes at a very small cost since it involves one solution of a KKT system of linear equations.
\end{myrem}

\section{Convergence results} {\label{sec:conv}}
In the previous section we gave an implementable algorithmic scheme that computes the solution to the CLQR problem. 
Here we provide all the necessary proofs which allow us to assert that the solution generated by Algorithm~\ref{al:AFBS_CLQR} via~(\ref{eq:primalSequence}) indeed converges to the true optimizer of the CLQR problem. In what follows $\lambdainf^\infty$ denotes any optimal solution to the dual problem~(\ref{eq:generic_dual}) (which exists under Assumption~\ref{as:basic_reg} but may not be unique) and $(\uinf^\infty, \xinf^\infty)$ the optimal solution to the primal problem~(\ref{mpc-initial_reg}). Our main result is:

\begin{thm}[Main Theorem]{\label{theor:master}}
Suppose Assumption~\ref{as:basic_reg} holds and let $\lambdainf^k$ be a sequence of iterates generated by Algorithm~\ref{al:AFBS} and $(\xinf^k,\uinf^k)$ the associated primal sequence given by~(\ref{eq:primalSequence}) and let $L$ be a Lipschitz constant of $\nabla h^\star(\cdot)$. The following statements hold:
\begin{enumerate}[label=(\roman{*}), ref=(\roman{*})]
\item The composite function $F(\lambdainf)=h^\star(\lambdainf)+\delta_-(\lambdainf)$ as defined in~(\ref{eq:generic_dual}) converges as
\[
	F(\lambdainf^k)-F(\lambdainf^\infty) \le \frac{a^2L}{2(k+a-1)^2}\|\lambdainf^0-\lambdainf^\infty\|_{\Hl}^2\enspace.
\]
\item The sequence of the dual iterates $\left(\lambdainf^k\right)_{k\in\natur}$ converges weakly (see Definition~\ref{def:weakconv1} in Appendix~\ref{App.A}) to an optimizer, that is,
\[
 \lambdainf^k\rightharpoonup \lambdainf^{\infty} 
\]
for some $\lambdainf^{\infty}\in\argmin(F)$.
\item The input sequence $\left(\uinf^k\right)_{k\in\natur}$ converges strongly to the unique minimizer as
\begin{equation*}
  \|\uinf^k-\uinf^\infty\|_{\Hu}\le a\sqrt{\frac{L}{\mu}}\frac{\|\lambdainf^0-\lambdainf^\infty\|_{\Hl}}{(k+a-1)}\enspace,
\end{equation*}
where $\mu>0$ is the strong convexity modulus of $f(\uinf)$.
\item The state sequence $\left(\xinf^k\right)_{k\in\natur}$ converges strongly to the unique minimizer as
\[
  \|\xinf^k-\xinf^\infty\|_{\Hx}\le a\sqrt{\frac{\|\Binf\|^2L}{\mu}}\frac{\|\lambdainf^0-\lambdainf^\infty\|_{\Hl}}{(k+a-1)}\enspace.
\]
\item The sequence $\left(T^k\right)_{k\in\natur}$ is bounded.
\end{enumerate}
\end{thm}
\begin{IEEEproof}
\begin{enumerate}[label=(\roman{*}), ref=(\roman{*})]
\item Convergence of $F(\lambdainf^k)$ with a constant stepsize is proven in~\cite[Theorem~1]{AFBS_chambolle}. Convergence at the same rate with an adaptive stepsize generated from the backtracking Algorithm~\ref{al:Backtracking} is proven in Lemma~\ref{lem::back_conv}, Appendix~\ref{App.D}.
\item The proof is stated in~\cite[Theorem~3]{AFBS_chambolle}.
\item The idea is to upper bound the input sequence's convergence rate making use of the result in point (ii). In order to do so we make use of strong duality. The proof is inspired from \cite[Theorem~4.1]{beck2014fast} and is as follows: \\
Let $\lambdainf^k\le 0$ generated from Step~7 of Algorithm~\ref{al:AFBS_CLQR}. Denote
\begin{equation}{\label{eq::fprime}}
\uinf^k = \underset{\uinf\in\Hu}{\argmin}\left\{f^\prime(\uinf):=f(\uinf)+\langle \lambdainf^k, \cinf-\Cinf\uinf\rangle\right\}\enspace,
\end{equation}
where $f(\uinf) = \frac{1}{2}\uinf^\star\Hinf\uinf  + \hinf^\star\uinf  + r$ as defined in Section~\ref{sec:intro}.
Then we have that the Lagrangian of~(\ref{eq:generic_pr}) evaluated at $\lambdainf^k$ is $\mathcal{L}(\uinf,\lambdainf^k)=f^\prime(\uinf)$.
The function $f(\uinf)$ is strongly convex with modulus $\mu\ge \lambda_{\mathrm{min}}(R)>0$, where $\lambda_{\mathrm{min}}(R)$ denotes the smallest eigenvalue of $R$. Strong convexity of $f^\prime(\uinf)$ with modulus $\mu$ follows directly. Using~(\ref{eq::fprime}), it holds that
\[
 f^\prime(\uinf) - f^\prime(\uinf^k) \ge \frac{\mu}{2}\|\uinf-\uinf^k\|_{\Hu}^2, \quad \forall \uinf\in\Hu\enspace,
\] 
or, equivalently,
\begin{equation}{\label{eq::Lagrange_strong_conv}}
 \mathcal{L}(\uinf,\lambdainf^k) - \mathcal{L}(\uinf^k,\lambdainf^k) \ge \frac{\mu}{2}\|\uinf-\uinf^k\|_{\Hu}^2, \quad \forall \uinf\in\Hu\enspace.
\end{equation}
Substituting $\uinf=\uinf^\infty$ in~(\ref{eq::Lagrange_strong_conv}) and by observing that $\underset{\lambdainf\le 0}{\max}\; \mathcal{L}(\uinf^\infty,\lambdainf)\ge \mathcal{L}(\uinf^\infty,\lambdainf^k)$, we have that
\begin{equation}{\label{eq::Lagrange_strong_conv2}}
  \mathcal{L}(\uinf^\infty,\lambdainf^\infty) - \mathcal{L}(\uinf^k,\lambdainf^k) \ge \frac{\mu}{2}\|\uinf^\infty-\uinf^k\|_{\Hu}^2, \quad \forall \uinf\in\Hu\enspace.
\end{equation}
We have managed to derive an upper bound for the distance of the generated sequence of primal minimizers $\left(\uinf^k\right)_{k\in\natur}$ from the optimal one. The last step is to show that the Lagrangian $\mathcal{L}(\uinf,\lambdainf)$ is associated to the composite objective $F(\lambdainf)$. This can be easily shown as follows:
\begin{align*}
 \mathcal{L}(\uinf^k,\lambdainf^k) &= \underset{\uinf\in\Hu}{\min}\left\{f(\uinf)+\langle \lambdainf^k, \cinf-\Cinf\uinf\rangle\right\} \\
 					     &= -\underset{\uinf\in\Hu}{\max}\left\{-f(\uinf)+\langle \lambdainf^k, \Cinf\uinf\rangle\right\} + \langle\lambdainf^k,\cinf\rangle \\
 					     &= -f^\star(\Cinf^\star \lambdainf^k) +  \langle \lambdainf^k, \cinf \rangle \\
 					     &= -F(\lambdainf^k), \quad \mbox{ by~(\ref{eq:generic_dual})}\enspace.
\end{align*} 
From strong duality and the fact that $-F(\lambdainf^k)$ converges to the optimal dual value (point (i)), we have that the optimal value of the dual function $-F(\lambdainf)$ coincides with that of the Lagrangian evaluated at the saddle point $(\uinf^\infty,\lambdainf^\infty)$, \ie, \mbox{$\mathcal{L}(\uinf^\infty,\lambdainf^\infty)=\underset{\lambdainf\in\Hl}{\max}\left\{-F(\lambdainf)\right\}=-F(\lambdainf^\infty)$} (see \cite[Section~5.5.5]{cvxbook}). Making use of point (i), inequality~(\ref{eq::Lagrange_strong_conv2}) becomes
\begin{equation}
\frac{\mu}{2}\|\uinf^k-\uinf^\infty\|_{\Hu}^2 \le F(\lambdainf^k) - F(\lambdainf^\infty) \le \frac{a^2L\|\lambdainf^0-\lambdainf^\infty\|_{\Hl}^2}{2(k+a-1)^2},
\end{equation}
which concludes the proof.
\item  The state sequence is generated by 
		 \begin{equation}{\label{eq::state_seq}}
		  \xinf^{k} = \Ainf x_{\mathrm{init}} + \Binf \uinf^{k}\enspace.
		 \end{equation}
		 Strong convergence of the input sequence $\left(\uinf^k\right)_{k\in\natur}$, along with the facts that $\Binf\colon\Hx\to\Hu$ is bounded (follows directly from Lemma~\ref{lem:bounded_C} in Appendix~C) and the uniqueness of $\uinf^\infty$ prove strong convergence of the state sequence with rate $1/k$, \ie,
		 \begin{align*}
		 \|\xinf^k-\xinf^\infty\|_{\Hx} &= \|\Binf (\uinf^{k}-\uinf^\infty)\|_{\Hx}  \\
		                          		&\le \|\Binf\|\|\uinf^{k}-\uinf^\infty\|_{\Hx} \\
		                         		&\le a\|\Binf\|\sqrt{\frac{L}{\mu}}\frac{\|\lambdainf^0-\lambdainf^\infty\|_{\Hl}}{(k+a-1)}\enspace,
		 \end{align*}
		with the last inequality following directly from point (iii).
\item 	The proof is already presented in~\cite{Sta_Kor_Jon_CLQR}, but we repeat it for completeness. First note that for the statement to hold it is 			sufficient to show that
		\begin{equation}\label{eq:limsupTbounded}
		\limsup_{k\to\infty} T^k < \infty.
		\end{equation}
		To prove~(\ref{eq:limsupTbounded}), define the sequence of the first hitting times of the interior of $\mathcal{S}$ as
		\[
		\tau^k := \inf\{i \ge T^k\mid x^k_i \in \interior \mathcal{S} \},\;\; k\in\mathbb{N}\cup\{+\infty\},
		\]
		where $\tau^\infty < \infty$ is the hitting time of the optimal state sequence~$\xinf^\infty$. Clearly, $\tau^k\ge T^k$ and $\tau^k < \infty$ since the origin is in the interior of $\mathcal{S}$ and for each $k\in\mathbb{N}$ the sequence $(\hat{x}^k_i)_{i\in\mathbb{N}}$ generated by Algorithm~\ref{al:AFBS_CLQR} converges to the origin as $i\to\infty$. We shall prove that $\limsup_{k\to\infty} \tau^k \le \tau^\infty < \infty$, which implies (\ref{eq:limsupTbounded}).

		The key piece for the proof is the weak convergence of $\hat{\xinf}^k$ to $\xinf^\infty$. This claim is proven in Lemma~\ref{lem:convx}, Appendix~\ref{App.B} and is a consequence of other intermediate results, all presented in the same Appendix.
		For the purpose of contradiction assume that there exists a subsequence $\tau^{k_j}$, $j\in\mathbb{N}$, with $\lim_{j\to\infty} \tau^{k_j} \ge \tau^\infty + 1$. Since the sequence of hitting times $\tau^k$ is integer valued, this implies that there exists a $j^\star\in \mathbb{N}$ such that $\tau^{k_j} \ge \tau^\infty + 1 $ for all $j\ge j^\star$. We now use this to contradict the weak convergence of $\hat{\xinf}^k$ to $\xinf^\infty$. To this end, observe that $x_{\tau^{\infty}}^\infty \in \interior\mathcal{S} $ whereas $\hat{x}_{\tau^{\infty}}^{k_j} \notin \interior\mathcal{S}$ for all $j\ge j^\star$. By the definition of the interior there exists an $\epsilon > 0$ such that $z\in \interior\mathcal{S}$ for all $z$ with $\|z-x^\infty_{\tau^\infty}\|_2 < \epsilon$. Therefore $\|\hat{x}_{\tau^{\infty}}^{k_j}-x^\infty_{\tau^\infty}\|_2 \ge \epsilon$ for all $j\ge j^\star$, and consequently 
		\begin{align*}
		\langle \hat{\xinf}^{k_j}-\xinf^\infty, \zinf \rangle & = \sum_{i=0}^{\infty} w^{i}(\hat{x}_i^{k_j}-x_i^\infty)^\tr z_i \\ & \ge w^{\tau^\infty} (\hat{x}_{\tau^\infty}^{k_j}-x_{\tau^\infty}^\infty)^\tr z_{\tau^\infty} \\ & \ge  w^{\tau^\infty}\epsilon^2 >0\enspace,
		\end{align*} 
		for a sequence $\zinf$ with $z_{\tau^\infty}=\hat{x}_{\tau^\infty}^{k_j}-x_{\tau^\infty}^\infty$ and for all $j\ge j^\star $, contradicting the weak convergence of $\hat{\xinf}^k$ to $\xinf^\infty$ asserted by Lemma~\ref{lem:convx}.
		\end{enumerate}
		\end{IEEEproof}

%\subsection{Nominal stability of the resulting controller}
%Finally, we prove stability of the closed-loop system under the control law resulting from Algorithm \ref{al:AMA_CLQR}.
%\begin{thm}
%  The origin is an exponentially stable equilibrium for the closed-loop system.
%\end{thm}
%\begin{IEEEproof}
%  The generated $\zinf^{\infty}(=\xinf^{\infty},\uinf^{\infty})$ sequence of Algorithm \ref{al:AMA_CLQR} respects both the dynamics and the inequality constraints, a fact that stems from the convergence of AMA. The optimal value of the objective is finite, implying that $\left(Q^{1/2}x^{\infty}_i\right)_{i\in\natur}\xrightarrow{i\to\infty} 0$ and $\left(R^{1/2}u_i^{\infty}\right)_{i\in\natur}\xrightarrow{i\to\infty} 0$ and from the positive definiteness of $Q$ and $R$ we have that $\left(x_i^{\infty}\right)_{i\in\natur}\xrightarrow{i\to\infty} 0$ and $\left(u_i^{\infty}\right)_{i\in\natur}\xrightarrow{i\to\infty} 0$. Since the nominal closed-loop trajectory is identical to the open-loop prediction, the closed-loop system is exponentially stable at the origin. 
%\end{IEEEproof}

\section{Computational aspects and warm-starting}{\label{sec:lin_sys_solve}}
Having presented the algorithm and its convergence results, we now focus on the computational aspects that render the algorithm a practical option to alternatives such as MPC. We start with explaining how no prior knowledge of a fixed stepsize is needed, give some references concerning the solution of the linear system, which is the most expensive operation of the method, and we conclude the section by suggesting a warm-starting scheme. 
\subsection{Stepsize selection}\label{sec:stepsize}
 The stepsize used in Algorithm~\ref{al:AFBS_CLQR} is computed as the reciprocal of the Lipschitz constant of $h^\star$. For the problem discussed here this can be explicitly computed. We have that 
 \begin{align}
 \|\nabla h^{\star}(\lambdainf_1)-\nabla h^{\star}(\lambdainf_2)\| \nonumber &= \|\Cinf\Hinf^{-1}\Cinf^{\star}\Winf(\lambdainf_1-\lambdainf_2)\| \\
                                                                   \nonumber &\leq \|\Cinf\Hinf^{-1}\Cinf^{\star}\Winf\|\|\lambdainf_1-\lambdainf_2\| \\
                                                                   \nonumber &= \|\Hinf^{-1}\|\|\Cinf\|^2\|\|\lambdainf_1-\lambdainf_2\|\enspace,
 \end{align}
 where $\Winf$ is the diagonal operator constituted of the decaying weighting sequence, \ie, $\Winf = \diag(I, wI, w^2 I,\ldots)$. The last equality follows from the fact that $\Winf$ contains a non-increasing sequence with the largest element being one. Hence $L(h^\star) = \|\Hinf^{-1}\|\|\Cinf\|^2$, which requires computation of the operator norm $\|\Hinf^{-1}\|$ and $\Cinf$. The proofs for boundedness of the two operators, as well as the computations of their bounds are derived in Appendix~\ref{App.C}, Lemmas~\ref{lem:bounded_C}~and~\ref{lem:bounded_Hinv}. 

 Although valid, this offline computation of the stepsize tends to be conservative in many cases, due to the conservativeness of the computed upper bounds. An elegant and practical method to achieve faster convergence is to employ an algorithm that locally estimates the Lipschitz constant online, at every iteration of Algorithm~\ref{al:AFBS_CLQR}. In order to do so, we use the \emph{backtracking stepsize rule} suggested in~\cite{fista}. The idea is simple: after each iteration of the algorithm, we make a quadratic approximation model of the function around the successor point, making use of the knowledge of the exact point and its gradient value. A quadratic term with varying curvature is added on top of the linear (first order Taylor) approximation, and the curvature is adapted recursively until our quadratic approximant upper bounds the original function, centered around the given point. Thus, the quadratic model's curvature is an estimate of the Lipschitz constant of the gradient of the original function. It is proven in~\cite{fista} that the locally evaluated Lipschitz constant $L^k$ is related to the global one by $\beta L(h^\star)\le L^k\le \gamma L(h^\star)$, where $\beta=\frac{L^0}{L(h^\star)}$ and $\gamma > 1$, $L^0>0$ being an initial estimate. Consequently the rule allows for smaller $L$'s and hence larger stepsizes, \ie, faster convergence in practice.

\subsection{Complexity}
The most expensive operation of Algorithm~\ref{al:AFBS_CLQR} is the linear system solve in Step~3. There is a variety of ways to perform this step, \ie, solving the KKT system or perform the Riccati recursion when both states and inputs are considered, invert the dense Hessian when only the inputs are considered. In the first case, a sparse (permuted) $\mathrm{LDL^\tr}$ factorization can be performed with cost $T(n+m)^3$ flops, followed by a forward-backward solve at $T(n+m)^2$ flops. This approach has been followed in~\cite{soc}. A discussion on the KKT system solve and the Riccati recursion approach is contained in~\cite{frison_riccati}, where the corresponding complexities are analyzed and compared in detail.

In the case of the condensed formulation, the linear system solve can be efficiently perfomed by first applying a Cholesky factorization on $H$ (being a finite truncation of the $\Hinf$ operator in~(\ref{eq:generic_pr})), followed by a forward-backward substitution, (see \cite[Appendix~C]{cvxbook}). Although the condensed form of the optimal control problem that is used in the derivations is, generally, not advised for long horizons, recent advancements can render this approach very efficient~\cite{Axehill201237},~\cite{FrisonJ13}. More specifically, the two proposed algorithms that perform factorization and solve of the condensed system in~\cite{FrisonJ13} come with a reduced complexity of $O(Tn^3)$ and $O(T^2n^2)$, respectively. 

Whether considering the sparse or the dense formulation, note that the factorization steps would have to be performed several times until the `correct' horizon $T^\infty$ has been identified, since the size of the corresponding matrices (KKT or Hessian $H$) increase as the algorithm progresses. This typically happens within the first few tens of iterations. An alternative would be to directly perform a matrix inversion (for small to medium sized systems) and only factorize once, when the horizon seems to have stabilised.

The backtracking scheme contributes to the complexity of Algorithm~\ref{al:AFBS_CLQR} by requiring a number of function evaluations per iteration, both for the quadratic model and for the original smooth function $h^\star(\cdot)$ (see~\cite{fista} for more details). The rest of the steps are simple vector updates of negligible cost.

\subsection{Warm-starting}
In the nominal case, \ie, when no noise and no model uncertainty are present, the open loop infinite-horizon control sequence generated from Algorithm~\ref{al:AFBS_CLQR} coincides with the control sequence generated by the optimal closed-loop feedback controller. Consequently, there is no need to re-optimize in a receding horizon fashion. Solving the CLQR problem just once is sufficient.

In the more realistic scenario where the predicted initial state differs from the measured one, the algorithm has to be re-applied. Provided that the prediction is not very different from the actual state, a good strategy is to initialize the decision variables (states and inputs) of the new problem with the values predicted from the previous one. This is commonly known as warm-starting. In our case, warm-starting has to be performed in the dual variables.

Note that once Algorithm~\ref{al:AFBS_CLQR} has run once, a hitting time $T^\infty$ has been generated, along with an optimal dual sequence $\lambda^\infty$ of corresponding length. It is expected that when computing the control law the hitting time should decrease by one at each solve, provided warm-starting from the optimal dual sequence that was generated once in the beginning. Hence, we suggest a heuristic scheme where the `constrained' (nonzero) part of the preceding shifted dual sequence is used to initialize each subsequent CLQR problem. The computation time thus reduces significantly, with the horizon practically shrinking to zero once the intial state is identified to be inside the maximal positively invariant set of the LQ controller.

Application of the scheme is presented in Section~\ref{sec:examples}. It is observed that it behaves particularly well for small perturbations of the initial state.

\section{Examples}{\label{sec:examples}}
For illustrative purposes, we run the algorithm on two systems, a small unstable system with two states and one input and a linearized model of a quadrocopter with 12 states and 4 inputs. We use the small example as a benchmark for graphical illustrations, while the larger one exhibits the computational efficiency of the proposed scheme. The comparison is performed against the same implementation of the AFBS algorithm for finite horizon lengths. It is of course understood that there exist several methods capable of solving a finite horizon MPC problems (see interior point, active set, etc.), among which, optimal first order methods have gained considerable attention over the last few years, rendering them a competitive alternative~\cite{fast_gradient},~\cite{ECC_14_short_survey}. Consequently, comparing against an optimal first order method provides a valid basis for evaluating the potential of our scheme. In the two examples the termination criterion is simply set as $\|\lambdainf^k-\hat{\lambdainf}^{k+1}\| \leq 10^{-4}$.
\subsection{Toy system}
Consider the following system defined as 
\[
A = \begin{bmatrix*}[r] 1.1&2\\0&0.95 \end{bmatrix*}, \quad B = \begin{bmatrix*}[r]0\\0.0787 \end{bmatrix*},
\]
\[
x_{i+1} = Ax_i + Bu_i, 
\]
with constraints 
\[
\|x\|_{\infty} \leq 10, \quad \|u\|_{\infty} \leq 1
\]
and $Q=\begin{bmatrix*}[r]2&-2\\-2&2 \end{bmatrix*}$, $R=2I$. Note that the system is unstable, hence a nonzero sequence of weights has to be chosen in order to ensure boundedness of the $\Cinf$ operator. We choose $w=1/1.1^2$ and the value of $a$ in Step~1 of Algorithm~\ref{al:AFBS_CLQR} is set to $4$.

\begin{figure}[!Htb]
      \begin{center}
      \includegraphics[scale=0.45]{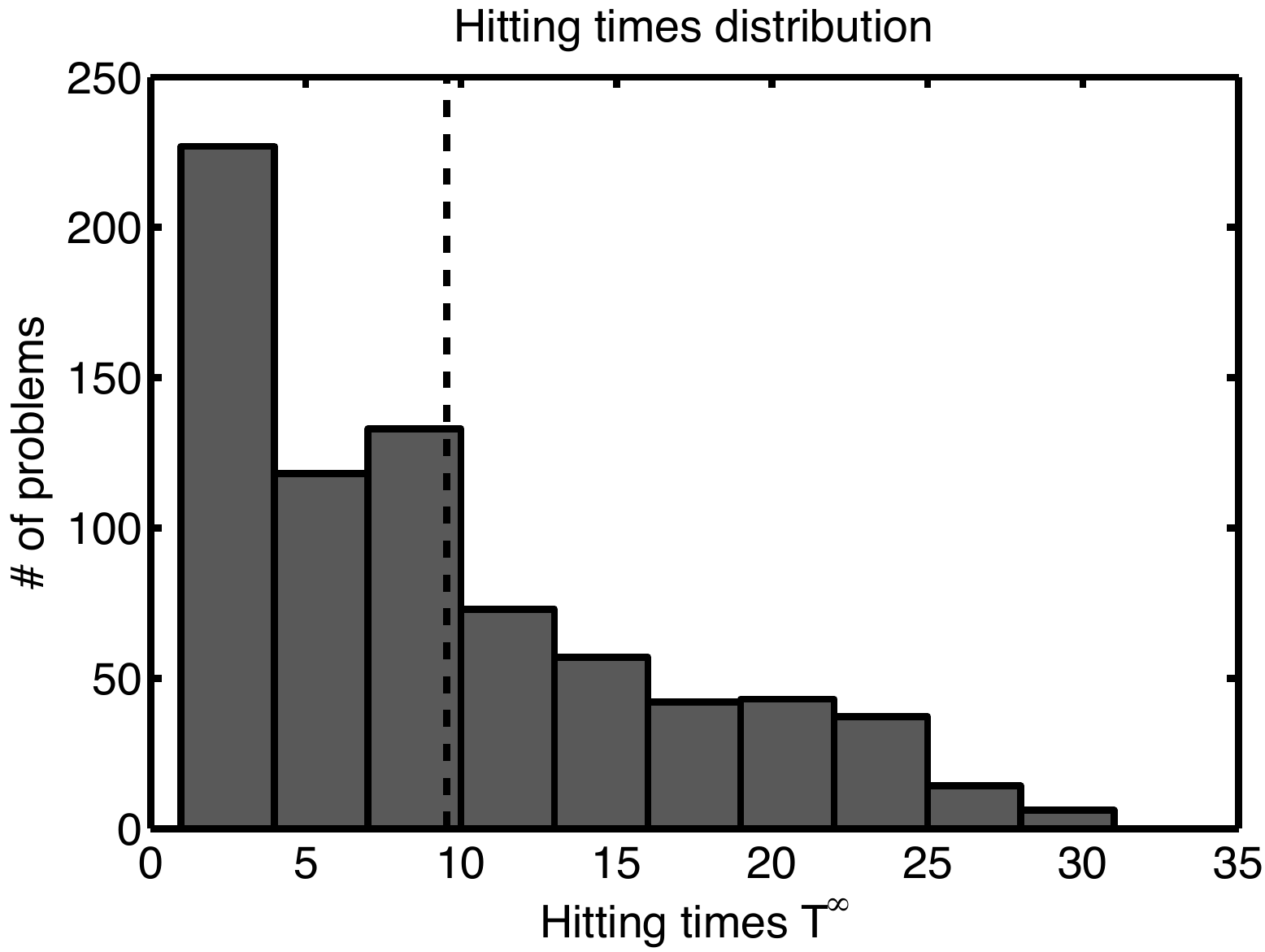}
      \caption{\small Histogram of $T^{\infty} = \max_{k}\{T^k\}$ for 750 initial conditions of the 2 state system sampled from a normal distribution around $(-3,0.3)$ with covariance matrix $\diag(4,0.4)$.}
      \label{fig:HTs_2states}
      \end{center}
  \end{figure}

The system is simulated for 750 different initial conditions $x_0$. In Figure~\ref{fig:HTs_2states} the distribution of $T^{\infty} = \max_k\{T^k\}$ is depicted. We see that $T^\infty$ goes up to 30, while the mean value is 9. 

We compare our method to the MPC approach both from the control and the algorithmic performance perspective. Regarding the former, we perform a comparison of the feasible sets for finite horizon implementations versus the maximal control invariant set in the case of CLQR, as well as the optimal value of the cost function. Regarding the latter, we perform comparisons in terms of the average number of iterations needed for convergence for several horizon length values in MPC versus the CLQR case. We also evaluate the conservativeness of our approach by computing the actual `optimal' horizon length for each initial condition we simulate. For both CLQR and MPC we make use of the same AFBS algorithm with backtracking employed and termination tolerance set to $10^{-4}$ as stated before.

\begin{figure}[!Htb]
      \begin{center}
      \includegraphics[scale=0.45]{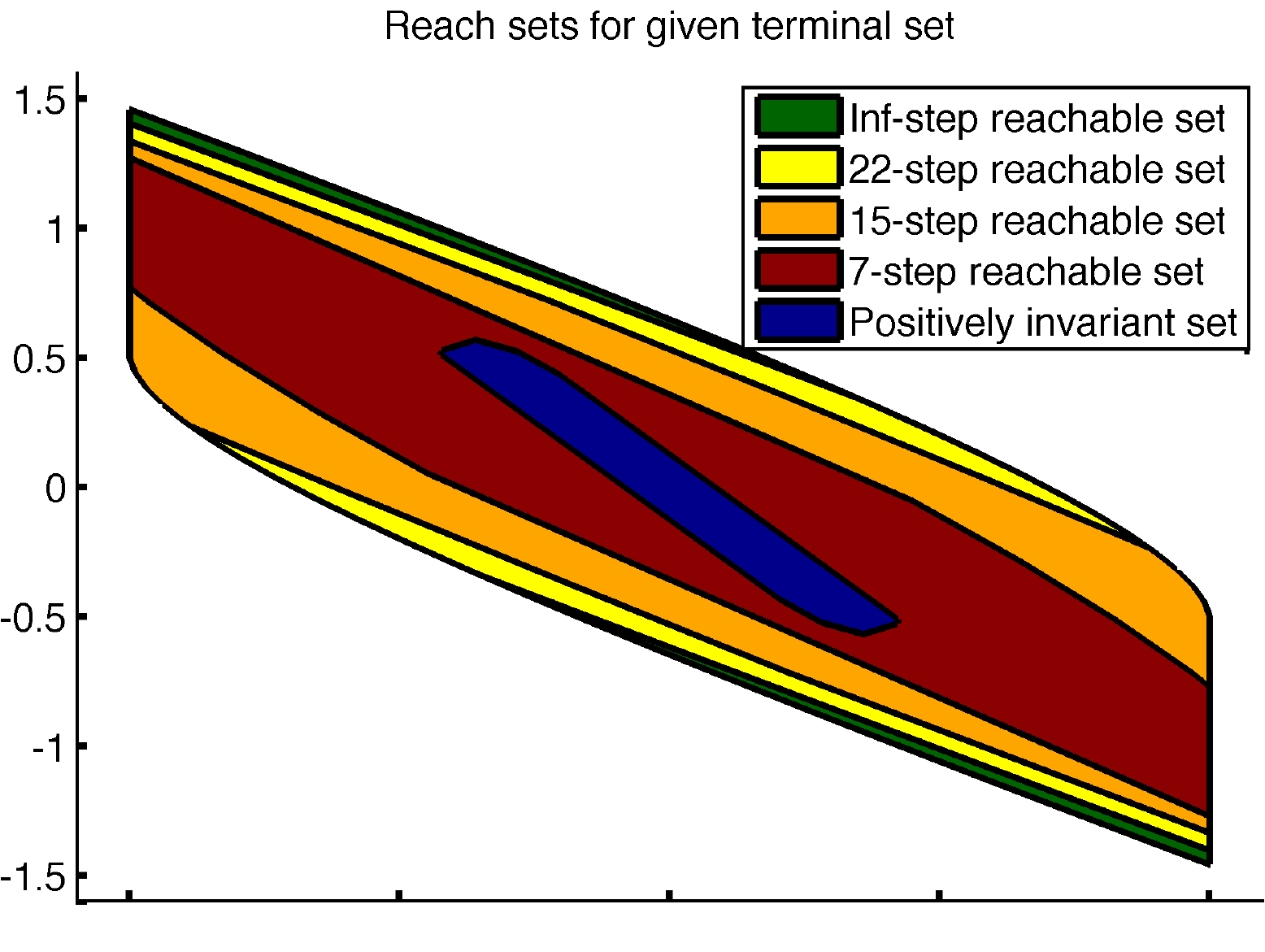}
      \caption{\small Reachable sets for several horizon lengths and the LQ terminal set. The computation was done using the MPT3 toolbox~\cite{MPT3}. It is apparent that a short horizon length reduces significantly the feasible region of the problem.}
      \label{fig:reach_sets}
      \end{center}
  \end{figure}

We perform the following simulations: we sample 243 initial conditions $x_{i,0}, \; i=1,\ldots,243$ from the maximal control invariant set (31-step) of the aforementioned constrained system (see Figure~\ref{fig:reach_sets}). For each point we compute:
\begin{enumerate}
\item The minimum horizon length $T_{\min}$, such that $x^\ast_{T_{\min}}(x_{i,0})$ resides in the maximum positively invariant set of the autonomous system $x^+ = (A+BK_{\mathrm{LQ}})x$, used as a terminal set.
\item The hitting time $T=T^\infty$ generated by our proposed scheme.
\end{enumerate}

\begin{figure}[!Htb]
      \begin{center}
      \includegraphics[scale=0.45]{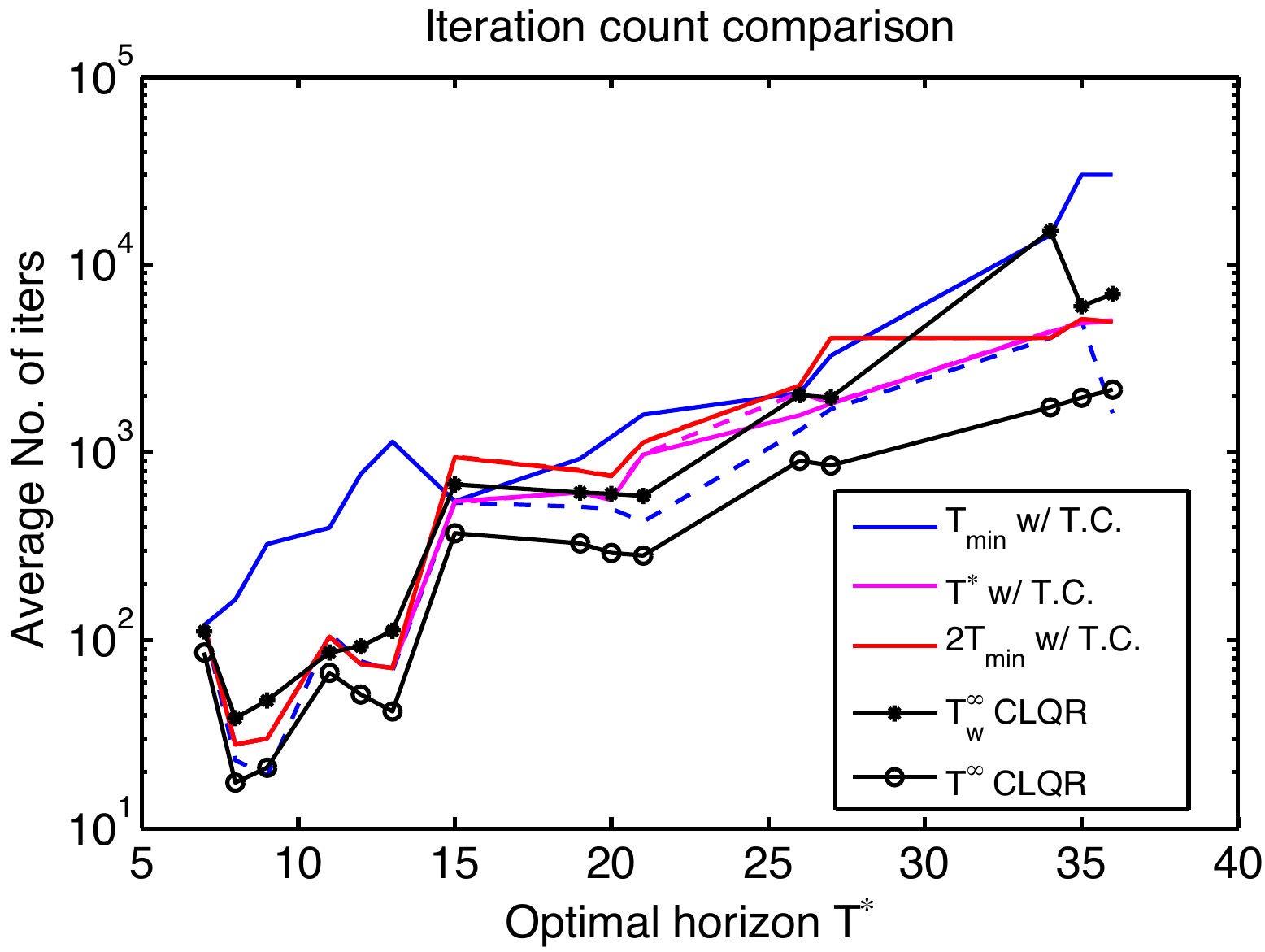}
      \caption{\small Comparison of MPC with finite horizon length, for several horizons, with CLQR. The x-axis corresponds to the optimal horizon length $T^\star$ per initial condition. As $T^\star$ increases, the corresponding states are sampled further from the origin. MPC with terminal set is depicted with the solid lines, while without terminal set with dashed lines. CLQR is performed both with the weight sequence, denoted as $T_w^\infty$ CLQR (provably convergent version), and without the weights, denoted as $T^\infty$ CLQR.}
      \label{fig:iters_MPC_CLQR}
      \end{center}
  \end{figure}

The following scenarios are generated: Firstly, an MPC problem with terminal set and horizon length $T_{\min}$ is solved, as described above. Subsequently we remove the terminal constraint and solve the MPC problem again for the same horizon length $T_{\min}$. We repeat the procedure described above (MPC with and without terminal set) for horizons $T=2T_{\min}$ and $T=T^\star$, where $T^\star$ is specified by identifying the first horizon length for which the feasible solution of the MPC problem had the terminal constraint inactive. Finally, a comparison with the proposed CLQR approach is performed, with horizon $T=T^\infty$. In this way, six MPC problems with finite horizon as well as the CLQR problem are solved, for each of the 243 initial conditions. The total number of iterations is averaged by the number of corresponding initial conditions with the same minimum horizon length $T_{\min}$. The results are presented in Figure~\ref{fig:iters_MPC_CLQR}.

The first observation from the plot is that inclusion of the terminal constraint generally increases the number of required iterations since more constraints become active at optimality. This is especially the case when the horizon is relatively short and the terminal constraint satisfaction is imposed as we see in the blue curves. As the horizon increases (red and magenta color), the terminal constraint might be inactive and there is no significant difference in the number of iterations in the two cases. 

A trend of increase in the number of iterations in all methods as $T^\star$ increases can be observed. This is expected since increase of $T^\star$ amounts to sampling of the initial state from regions of the feasible set that are further away from the origin. A consequence of the latter is that more constraints become active at optimality. 

Figure~\ref{fig:iters_MPC_CLQR} depicts a generally expected behavior. Short horizons in combination with a terminal set increase significantly the iteration count. As the horizon increases the iteration count drops. However, an unnecessarily long horizon (\eg, $2T_{\min}$) leads as well to an increase. In that sense, one can observe that MPC with $T=T^\star$ behaves better than the other lengths. The curve corresponding to CLQR is denoted with black lines with asterisks. We observe that the method is comparable to the MPC approaches where $T=T_{\min}$ and $T=T^\star$ with terminal set, performing relatively better in small to medium horizon lengths and worse for larger horizon lengths. The reason for the latter is that the weighting operator $\Winf$ becomes quickly ill-conditioned as $T$ grows, since its diagonal elements decay exponentially. This leads to ill-conditioning of the finite-dimensional truncation of the dual problem~(\ref{eq:generic_dual}) since $\Winf$ preconditions (through the weighted $l_2$ inner product~(\ref{def::l2_space})) the objective of~(\ref{eq:generic_dual}). As is commonly known first-order methods struggle with ill-conditioning. Efficient preconditioning of the operators should improve this and is a topic of further investigation.

We further perform the following simulation: The CLQR approach is applied, this time with the weights being set to one. The result is depicted as a black line with circles. In this case, CLQR clearly outperforms all MPC approaches with terminal set, as well as most of those instances without terminal set. Further simulations suggest that the proposed method performs very well, were we to drop the weights or when dealing with stable systems, when no scaling has to be performed. 

\begin{figure}[!Htb]
      \begin{center}
      \includegraphics[scale=0.45]{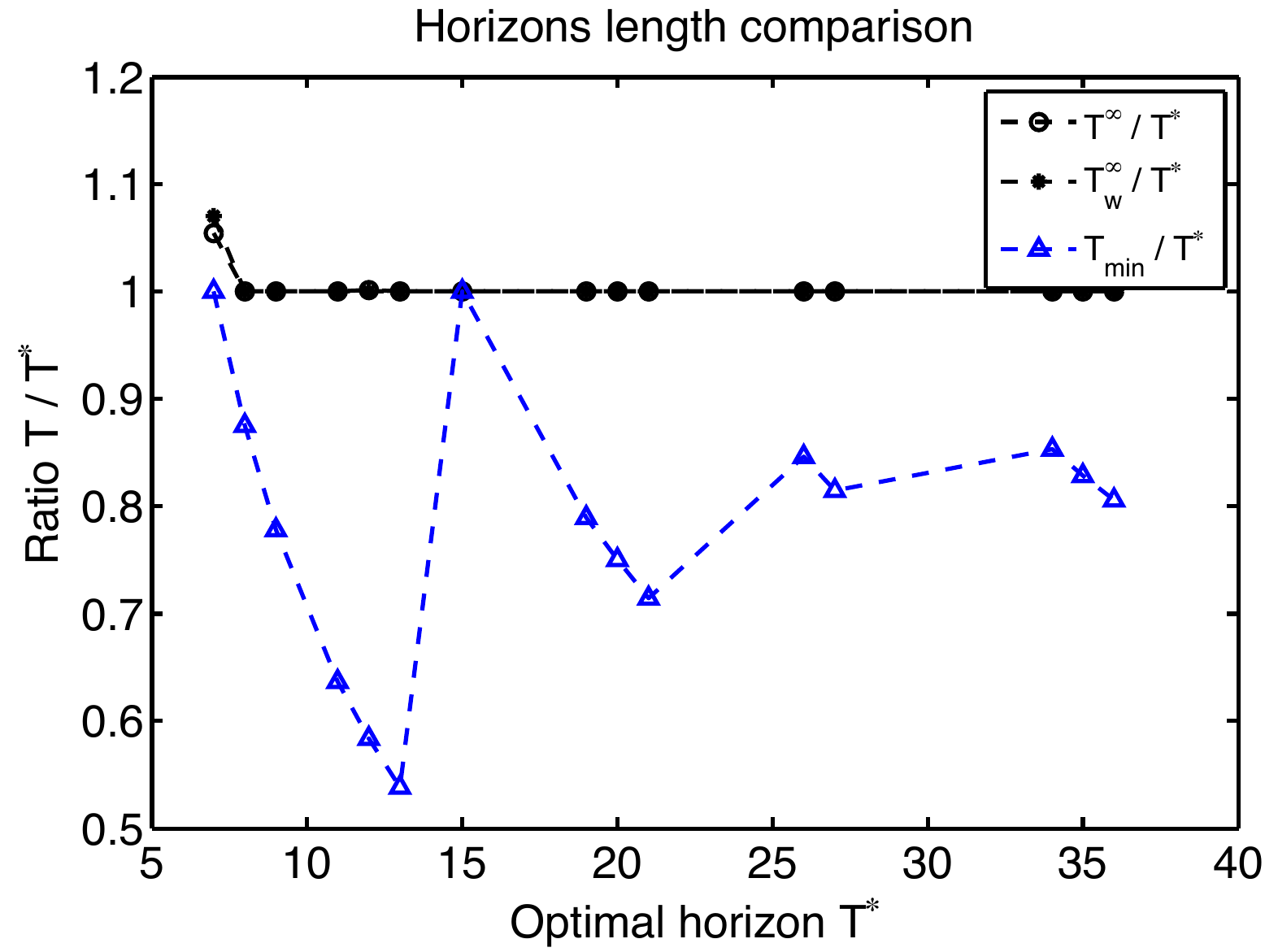}
      \caption{\small Evolution of the ratios $\frac{T_{\min}}{T^\star}$, $\frac{T_w^\infty}{T^\star}$ and $\frac{T^\infty}{T^\star}$ for the sampled initial conditions.}
      \label{fig:horzs_MPC_CLQR}
      \end{center}
  \end{figure}

Another interesting point is that the average hitting times $T^\infty$ generated from Algorithm~\ref{al:AFBS_CLQR} are almost identical to the averaged optimal ones $T^\star$, with a very slight increase. This fact is depicted in Figure~\ref{fig:horzs_MPC_CLQR}, where the ratio is always very close to one. On the contrary, the minimum required horizon length $T_{\min}$ is observed to be up to $45\%$ smaller than the optimal length $T^\star$ in some cases, leading to a significant increase in the objective's cost when compared to the infinite horizon approach.

\iffalse
	\begin{table}
	\normalsize
	\begin{center}
	\scalebox{0.75}{
	\begin{tabular}{|c||c|c|c|c|}
	\hline
	 & $\mathcal{S}_{6,7}$ & $\mathcal{S}_{14,15}$ & $\mathcal{S}_{21,22}$ & $\mathcal{S}_{28,29}$ \\
	\hline\hline
	Av. Iters. AFBS for CLQR & 90.31 & 486.34  & 1238.4 & 2289.4 \\
	\hline
	Av. Iters. AFBS for MPC w/ T.C. & 248.40 & 715.66 & 2634.7 & 5726.3 \\
	\hline
	Av. Iters. AFBS for MPC w/o T.C. & 95.77 & 566.09 & 1847.2 & 4309.6 \\
	\hline
	Av. Hor. Lenth $T^\infty$  & 10.06 & 18.18 & 27.20 & 34.75\\
	\hline
	Av. Opt. Hor. Lenth $T^\star$  & 9.99 & 18.06 & 26.84 & 24.75\\
	\hline
	Av. Optimal Val. AFBS for CLQR & 292.69 & 599.66  & 1079.4 & 1668.3 \\
	\hline
	Av. Optimal Val. AFBS for MPC w/ T.C. & 292.72 & 599.66 & 1079.4 & 1668.3 \\
	\hline
	\end{tabular}
	}
	\end{center}
	\caption{Comparison of Algorithm~\ref{al:AFBS_CLQR} (CLQR) to MPC solution for the two-states system.}
	\label{TableI}
	\end{table}
\fi

  \subsection{Quadcopter system}
  The next system we consider is a quadcopter linearized in a hovering equilibrium. The system has 12 states which correspond to position, angle and the corresponding velocities. There are four inputs corresponding to the four propellers.
  There are box constraints on all states and inputs, mainly ensuring the validity of the linearized model. The system is marginally stable; thus the weight $w$ is set to one.

\begin{figure}[!Htb]
      \begin{center}
      \includegraphics[scale=0.45]{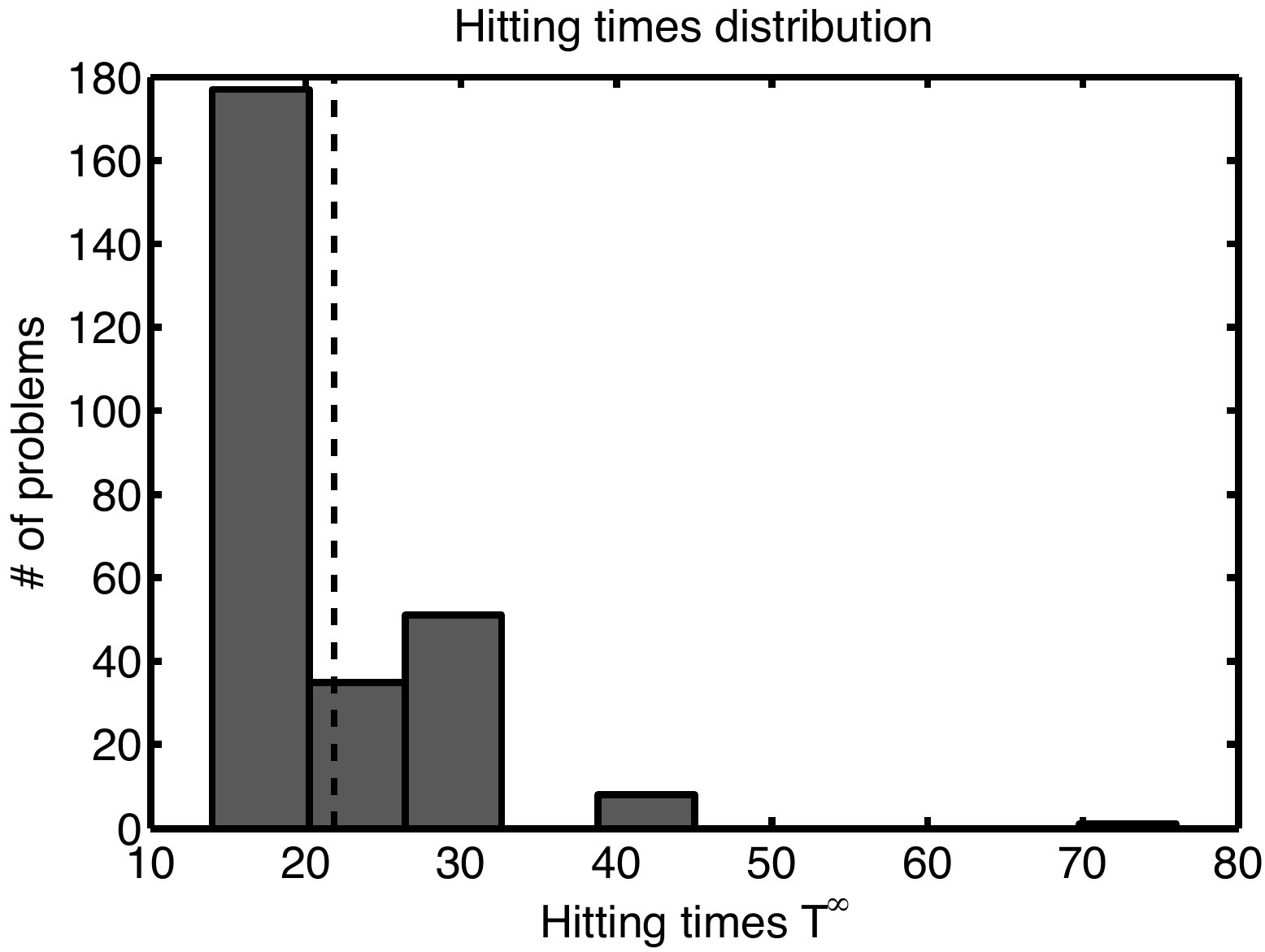}
      \caption{\small Histogram of $T^{\infty} = \max_{k}\{T^k\}$ for 272 initial conditions sampled with a Hit-And-Run algorithm.}
      \label{fig:HTs_quad}
      \end{center}
  \end{figure}

 We simulate the algorithm starting from initial conditions randomly selected as follows: starting from a random feasible initial condition, we generate random directions on a unit ball centered around it and sample points along each of them. The points are generated from a normal distribution with standard deviation 0.15. Finally, we keep the initial conditions that result in feasible closed loop problems. The result of this step is 272 feasible initial conditions for the CLQR Algorithm~\ref{al:AFBS_CLQR}. A histogram of $T^\infty = \max_k\{T^k\}$ is presented in Figure \ref{fig:HTs_quad}.  

\begin{figure}[!Htb]
      \begin{center}
      \includegraphics[scale=0.45]{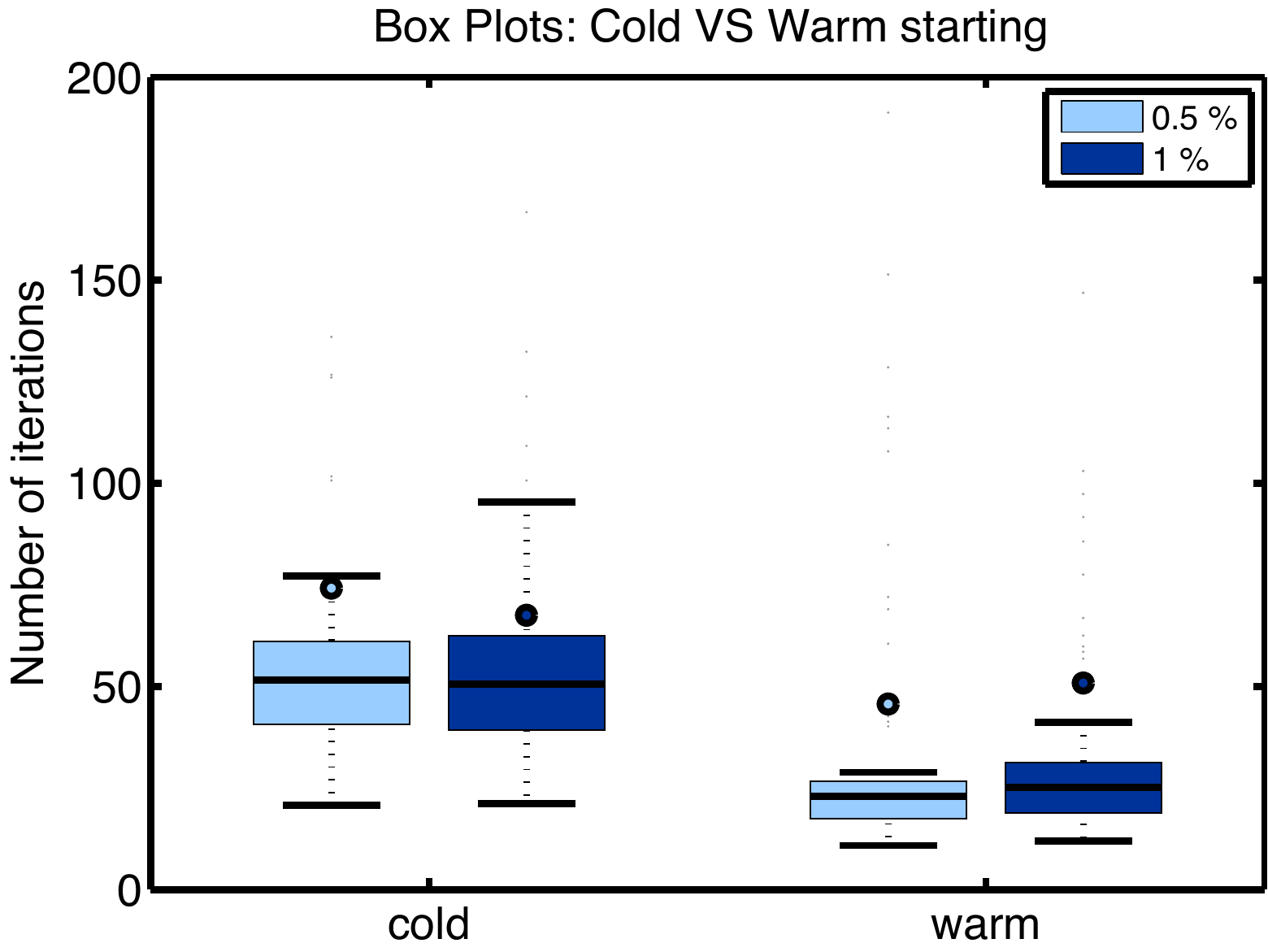}
      \caption{\small `Box $\&$ Whiskers' plot for the average number of iterations of the warm-starting policy, given $0.5\%$ and $1\%$ uniform perturbations of the initial state. The horizontal line inside the box corresponds to the median, the edges of the box are the $25^{\mathrm{th}}$ and $75^{\mathrm{th}}$ percentiles while the horizontal lines outside the boxes correspond to the most extreme data points not considered outliers. Finally, the colored dots correspond to mean values. One can observe that in the case of warm-starting the means are located outside the extreme points band, suggesting the existence of large outliers. However, warm-starting improves the performance of the suggested method in both cases, in all the depicted statistical measures.}
      \label{fig::box_whiskers}
      \end{center}
  \end{figure}

We conclude the section by applying the warm-starting heuristic scheme suggested in Section~\ref{sec:lin_sys_solve}. We consider two scenarios; in the first one we uniformly perturb the initial state by $0.5\%$ of its nominal value and run Algorithm~\ref{al:AFBS_CLQR} in closed loop for 78 different initial conditions. In the second scenario we perturb the state by $1\%$, for 68 different initial conditions. We solve 15 consecutive problems per initial condition and subsequently compare the average number of iterations as well as the generated hitting times per problem solve with and without warm-starting the dual variables. The results are summarized in Figure~\ref{fig::box_whiskers}. It is evident that warm-strating consistently reduces the number of iterations in both cases.

\begin{figure}[!Htb]
      \begin{center}
      \includegraphics[scale=0.45]{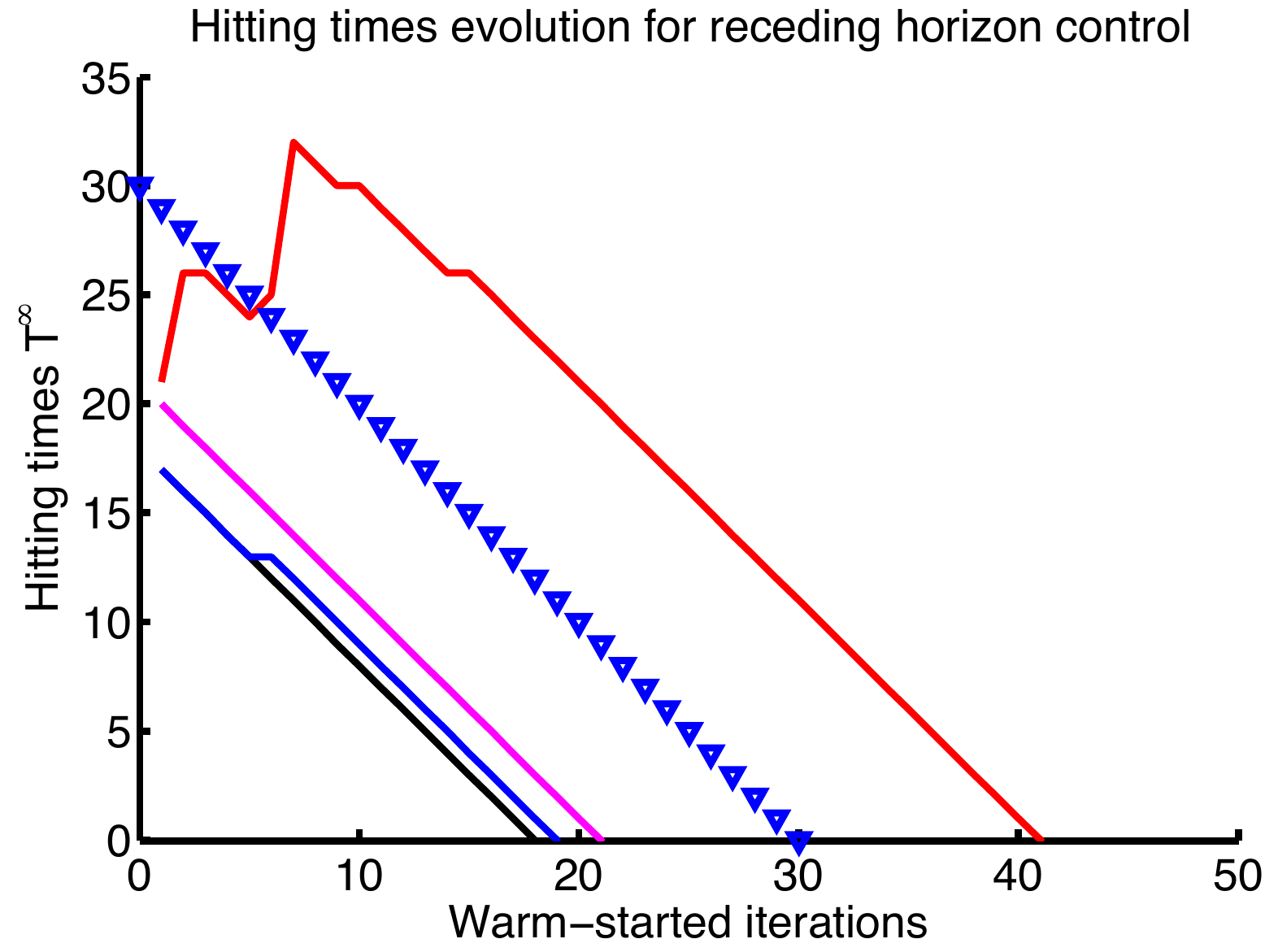}
      \caption{\small Hitting times evolution for 50 subsequent solves, warm-started at the previous (shifted) dual optimizer. Four different initial conditions are depicted with the solid lines. The linear rate is depicted (for comparison) with the blue triangles. Note that the horizon decreases almost linearly, resulting from the relatively small perturbations of the initial state. The case depicted in red corresponds to an outlier where the perturbed initial condition rendered the previously computed horizon an infeasible option, leading to an increase of its length.}
      \label{fig:warm_horizon}
      \end{center}
  \end{figure}

  Finally, the active reduction in the horizon length thanks to warm-starting is illustrated in Figure~\ref{fig:warm_horizon} for four different initial states in the case of the $1\%$ perturbation. There is an (almost monotonic) decreasing trend, with the horizon finally shrinking to zero when the initial state resides in the positively invariant set.

\section{Conclusion}
This work presented an algorithmic scheme capable of solving the constrained linear quadratic regulator problem in real time. The algorithm is an accelerated version of Forward-Backward Splitting, a popular proximal gradient scheme. The  approach is to write the problem in its condensed form and dualize, which leads to the minimization of an infinite-dimensional quadratic functional subject to non-positivity constraints. The resulting infinite dimensional problem can be tackled in finite dimensions by observing that the dual sequence has always only finitely many non-zero elements. The proposed algorithm makes no use of terminal invariant sets and provably converges to the optimal solution of the infinite-horizon problem. Regarding the implementation aspects, the algorithm can be highly competitive since it enjoys the convergence properties of optimal first order methods whose each iteration is computationally cheap. In addition, it requires minimal a priori information since the most crucial quantities are computed online, and there are no unreasonable or conservative assumptions on the problem's structure. On the other hand, the number of iterations needed for convergence can be quite large for very unstable (large) systems. To this end, preconditioning of the problem has to be considered, with the obvious difficulty of computing preconditioners for infinite dimensional operators. Some steps towards this direction are taken in, \eg, \cite{vmfbs}, but this problem remains a topic of future research. 

The algorithm used, AFBS, requires the objective function to be written as the sum of a smooth and a (possibly) nonsmooth function whose proximal operator is inexpensive to compute. Smoothness of the dual problem is recovered from strong convexity of the primal. The polytopic constraint set becomes a nonnegativity constraint when looking at the dual variables, resulting in a simple proximal operator. %Furthermore, if the constraints imposed on state and inputs are simple (\eg, boxes), there is no need to dualize. One can directly tackle the problem in its primal variables, with the apparent advantage of having less variables to optimize over.

A direct extension of the existing scheme is the case of (some states) tracking a reference signal. By casting the tracking problem in a delta-formulation with respect to the inputs and states, the resulting problem retains the same structure as the regulation problem, with the difference of varying constraint sets due to the generated steady-state points.

\appendices
\section{Required Operator Theory}{\label{App.A}}
For the sake of clarity we introduce some notation and definitions from operator theory. The subsequent results hold for general real Hilbert spaces, including the special case of $\lw$ we consider. For an in-depth treatment of (monotone) operator theory, the interested reader is referred to \cite{book_comb} and \cite{Eckstein}.
We write variables in normal font, and we use the bold font to describe the infinite-dimensional variables we are manipulating in our problem description. 

\begin{mydef}{\label{def:wl2_hilbert_space}}
The $l^2$-weighted (or $\lw$) real Hilbert space $\mathcal{H}$ is defined by
\[
 \mathcal{H} = \left\{\zinf=\left(z_i\right)_{i \in \natur}\;:\;\sum_{i=0}^\infty\|z_i\|^2_2w^i<\infty\right\}, \quad w>0\enspace.
\]
\end{mydef}

%An \emph{operator} $A\colon \HI\to 2^{\HII}$ is a point-to-set map, \ie, $A$ maps every point $x\in\HI$ to a set $A(x)\subseteq\HII$. The operator is characterized by its \emph{graph}, $\gra A=\left\{(x,u)\in \HI\times\HII \mid u \in A(x)\right\}$. The \emph{inverse} $A^{-1}$ of $A$ is defined through its graph as  \\ $\gra A^{-1}=\left\{(u,x)\in\HII\times\HI\mid (x,u)\in\gra A\right\}$. Composition, scalar multiplication and addition of operators are well-defined operations.

\begin{mydef}{\label{def:boundedoper}}
A linear operator (mapping) $F\colon \HI\rightarrow \HII$ between two Hilbert spaces is said to be bounded if the operator norm $\|F\|$ of $F$, defined as
\[
\|F\|\coloneqq \sup_{\|x\|_{\HI}=1}\|Fx\|_{\HII}\enspace,
\]
satisfies $\|F\|<\infty$. The set of bounded operators between two Hilbert spaces $\HI$ and $\HII$ is denoted as $\mathcal{B}(\HI,\HII)$. 
\end{mydef}

%\begin{thm}[Bounded Inverse Theorem]{\label{def:boundedinverse}}
% Let $F\colon \HI\rightarrow \HII$ be a linear operator which is bijective, and $F\in\mathcal{B}(\HI,\HII)$. Then $F^{-1}\in\mathcal{B}(\HII,\HI)$.
%\end{thm}

\begin{thm}{\label{thm:adjoint}}
Let $\HI,\HII$ be real Hilbert spaces and $F\in\mathcal{B}(\HI,\HII)$. The adjoint of $F$ is the unique operator $F^{\star}\in\mathcal{B}(\HII,\HI)$ that satisfies
\[
 \langle Fx,y\rangle  = \langle x,F^{\star}y\rangle  \quad \forall x\in\HI, \; \forall y\in\HII \enspace.
\]
Moreover, $\|F\|=\|F^\star\|$.
\end{thm}

%\begin{prop}{\label{prop:subseq}}
%If $\left(\alpha^k\right)_{k\in J}$, where $J$ is a countable index set, is a convergent sequence, then every subsequence of that sequence converges to the same limit.
%\end{prop}

Subsequently, we introduce the notions of weak and strong convergence.
\begin{mydef}{\label{def:weakconv1}}
Let $\mathcal{H}$ be a Hilbert space. We say that $\left(x^k\right)_{k\in\natur}$ converges weakly to $x$ if $\forall y \in \mathcal{H}$ $\langle y,x^k\rangle \xrightarrow{k\to\infty}\langle y,x\rangle $. We denote weak convergence as $x^k \rightharpoonup x$.
\end{mydef}

\begin{mydef}{\label{def:strongconv}}
Let $(x^k)_{k\in\natur}$ be a sequence in $\mathcal{H}$. Then $\left(x^k\right)_{k\in\natur}$ converges strongly to $x$ if $\|x^k-x\|\xrightarrow{k\to\infty} 0$. We denote strong convergence as $x^k \rightarrow x$.
\end{mydef}

\begin{mydef}{\label{def:posdef}}
An operator $F: \mathcal{H}\to \mathcal{H}$ is \emph{positive-definite} if it is bounded, $F=F^{\star}$ and $\langle Fx,x \rangle \ge \alpha \|x \|_\mathcal{H} $ for some $\alpha>0$, for all $x \in \mathcal{H}$\enspace.
\end{mydef}
We recall that a positive definite operator $F$ is invertible and the inverse operator $F^{-1}$ is bounded.

\begin{mydef}{\label{def:proximal}} The proximity operator of a proper lower semi-continuous convex function $f\colon\mathcal{H}\rightarrow\reals$ is
\[
\prox{_f}(y) = \argmin_{x \in \mathcal{H}} \left\{f(x) + \frac{1}{2}\|x-y\|^2\right\}\enspace.
\]
\end{mydef}

\section{Convergence of the accelerated input-state sequences of Algorithm~\ref{al:AFBS_CLQR}}{\label{App.B}}
In this section we present the result that allows us to prove boundedness of the horizon sequence $T^k$, point (v) of Theorem~\ref{theor:master} in Section~\ref{sec:conv}. It is sufficient to show weak convergence of the $\hat \xinf^k$ sequence to the state optimizer $\xinf^\infty$, which is required in the proof of point (v) of Theorem~\ref{theor:master}.

The theorem is stated at the end of the Appendix as a sequence of intermediate results. We show that:
\begin{enumerate}
 \item The relaxed dual sequence $\hat{\lambdainf}^k$ converges weakly to a dual minimizer $\lambdainf^\infty$.
 \item Provided that the operator $\Cinf$ is bounded, the sequence $\Hinf\hat{\uinf}^k (= \Cinf^{\star}\hat{\lambdainf}^{k-1}-\hinf)$ converges weakly to $\Hinf\uinf^k$, and $\hat{\uinf}^k$ converges weakly to $\uinf^k$. From strong duality, $\uinf^k$ converges to the primal optimizer $\uinf^\infty$.
 \item Weak convergence of the accelerated state sequence $\hat{\xinf}^k$ to $\xinf^\infty$ follows directly.
 \end{enumerate}

Weak convergence of the relaxed sequence $\hat{\lambdainf}^k$ follows from Corollary~2 of~\cite{AFBS_chambolle}, which states that the error sequence $\left(\|\lambdainf^k-\lambdainf^{k-1}\|^2\right)_{k\in\natur}$ converges to zero with rate $1/k^2$. We state the result below.
 \begin{mylem}{\label{relaxed_dual_conv}}
  The relaxed sequence $\hat{\lambdainf}^k$ converges weakly to $\lambdainf^{\infty}$.
 \end{mylem}
 \begin{IEEEproof}
  Since $\|\lambdainf^k-\lambdainf^{k-1}\|^2 \to 0$ and $\alpha_k$ is bounded we also have $\nuinf^k = \sqrt{\alpha^k}(\lambdainf^k-\lambdainf^{k-1}) \to 0$. Since strong convergence implies weak convergence we have that $\langle \nuinf^k,\yinf\rangle \xrightarrow{k\to\infty}0 , \; \forall \yinf\in\Hl$.
  The relaxed sequence of duals $\hat{\lambdainf}$ can be written as $\hat{\lambdainf}^k = \lambdainf^k + \sqrt{\alpha^k}\nuinf^k$. Consequently, since $\lambdainf^k \rightharpoonup \lambdainf^\infty $, we have that 
  \begin{align*}
  \langle \hat{\lambdainf}^k,\yinf\rangle &=\langle \lambdainf^k + \sqrt{\alpha^k}\nuinf^k,\yinf\rangle \\
                                          &=\langle \lambdainf^k, \yinf\rangle + \sqrt{\alpha^k} \langle \nuinf^k,\yinf\rangle \to \lambdainf^\infty
  \end{align*}
for all $y \in \mathcal{H}_{\lambdainf}$ and hence $\hat\lambdainf^k \rightharpoonup \lambdainf^\infty$.
 \end{IEEEproof}

\begin{mylem}\label{lem:convu}
 The sequence $\left(\hat{\uinf}^k\right)_{k\in\natur}$ converges weakly to $\uinf^\infty$.
\end{mylem}
\begin{IEEEproof}
Writing down the relation between $\hat \uinf$ and $\hat\lambdainf$ from Lemma~\ref{lem:grad} in terms of the operators, we get $\hat \uinf^k = \Hinf^{-1}(\Cinf^\star \Winf \hat\lambdainf^k - \hinf)$. Similarly we have $ \uinf^k = \Hinf^{-1}(\Cinf^\star \Winf \lambdainf^k - \hinf)$. Now, from Theorem~\ref{theor:master} we have $\uinf^k \to \uinf^\infty$ and from Lemma~\ref{relaxed_dual_conv} we have $\hat\lambdainf^k - \lambdainf^k \rightharpoonup 0 $. Therefore, since $\Cinf^\star$, $\Winf$ and $\Hinf^{-1}$ are bounded operators (see Lemmas~\ref{lem:bounded_C} and~\ref{lem:bounded_Hinv} in Appendix~C) and since weak convergence is preserved under bounded linear mappings, we conclude that $\hat\uinf^k \rightharpoonup \uinf^\infty$.

\end{IEEEproof}

\begin{mylem}\label{lem:convx}
 The sequence $\left(\hat{\xinf}^k\right)_{k\in\natur}$ converges weakly to $\xinf^\infty$.
\end{mylem}
\begin{IEEEproof}
 Exactly as we did at point (iv) of Theorem~\ref{theor:master}, the accelerated state sequence can be written as
		 \[
		  \hat{\xinf}^{k} = \Ainf x_{\mathrm{init}} + \Binf \hat{\uinf}^{k}\enspace.
		 \]
		 Weak convergence $\hat \uinf^k \rightharpoonup \uinf^\infty $ and boundedness of $\Binf$ prove weak convergence of the accelerated state sequence to $\xinf^\infty$.
\end{IEEEproof}

\section{Boundedness of several operators}{\label{App.C}}
\begin{mylem}{\label{lem:bounded_C}}
 The operators $\Cinf$ and $\Winf$ are bounded.
\end{mylem}
\begin{IEEEproof}
Boundedness of $\Winf$ is trivial since it is a diagonal operator with non-increasing elements on the diagonal.

 The operator $\Cinf$ can be expressed as the following sum:
 \[
  \Cinf = \left[\begin{array}{ccc}
  C_u & 0 & \cdots \\ 0 & 0 & \cdots \\ 0 & C_u & \cdots \\ 0 & 0 & \cdots \\ \vdots&\vdots&\ddots\end{array}\right] + 
          \underbrace{\left[\begin{array}{ccc}
  0&0&\cdots \\ C_xB&0&\cdots \\ 0&0&\cdots \\ C_xA B&C_xB&\cdots\\ \vdots&\vdots&\ddots \end{array}\right]}_{\Cxinf}\enspace,
\]
 and by using the triangle inequality, we have
 \begin{equation}{\label{App.C:toprove}}
  \|\Cinf\| \le \sigma_{\max}(C_u) + \|\Cxinf\|.
 \end{equation}
 We thus have to show that $\yinf=\Cxinf\uinf$ is bounded, \ie,
 \[
 \underset{\|\uinf\|_{\Hu}= 1}{\sup}\|\yinf\|_{\Hl} < \infty.
 \]
 In order not to carry the zero rows of $\Cxinf$, we define $\bar{\yinf}\in\bar{\Hl}$ by dropping the zero elements
 of $\yinf$. This infinite-dimensional vector consists of the elements $\bar{y}_i=[\Cxinf]_i\uinf=\sum_{j=0}^{i-1} C_xA^{i-j-1}Bu_j\in\reals^{p_x}$. Note that 
\[
  \underset{\|\uinf\|_{\Hu}= 1}{\sup}\|\bar{\yinf}\|_{\bar{\Hl}} = \underset{\|\uinf\|_{\Hu}= 1}{\sup}\|\yinf\|_{\Hl}\enspace.
\]
 Focusing on the operator of interest, we have that
 \begin{equation}{\label{App.C::eq1}}
  \|\bar{\yinf}\|_{\bar{\Hl}} = \sqrt{\sum_{i=1}^\infty\|\bar{y}_i\|^2_2w^i} = \sqrt{\sum_{i=1}^\infty w^i\|\sum_{j=0}^{i-1} C_xA^{i-j-1}Bu_j\|^2_2}\enspace.
 \end{equation}
 Then
\begin{align*}
 \|\bar{\yinf}\|^2_{\bar{\Hl}} &= \sum_{i=1}^\infty w^i\|\sum_{j=0}^{i-1} C_xA^{i-j-1}Bu_j\|^2_2 \\
                         &= \sum_{i=1}^\infty \|\sum_{j=0}^{i-1} C_x(A\hat{w})^{i-j-1}Bu_j\hat{w}^{j+1}\|^2_2 \\ 
                         &= \hat{w}^2\sum_{i=1}^\infty \|\sum_{j=0}^{i-1} C_x\hat{A}^{i-j-1}B\hat{u}_j\|^2_2\enspace,
\end{align*}
where we introduced $\hat{w}=w^{1/2}$, $\hat{u}_j=u_j\hat{w}^j$ with $\hat{\uinf}\in l^2$ and $\hat{A}=A\hat{w}$.
Observing the above expression, one can identify that $\sum_{j=0}^{i-1} C_x\hat{A}^{i-j-1}B\hat{u}_j$ is the \emph{convolution sum} of the impulse response of the system 
\[
\Sigma \coloneqq \left(
\begin{array}{c|c}
\hat{A} & B \\ \hline
C_x & 0
\end{array}\right)
\]
with an input $\hat{u}$. More specifically, borrowing the notation from \cite{Antoulas} we denote the impulse response of $\Sigma$ as $h_{\Sigma,i} = C_x\hat{A}^{i-1}Bs_i$, where $s_i=1, \; i\geq 0$ is the unit step function. Then the convolution operator is defined
as the linear map 
$S_\Sigma\colon\hat{\uinf}\rightarrow \bar{\yinf}$ with $\bar{y}_i = \left(h_\Sigma\ast\hat{\uinf}\right)_i = \sum_{j=0}^{i-1}h_{\Sigma,i-j}\hat{u}_j = \left(S_\Sigma\hat{\uinf}\right)_i$. Thus we have that 
\begin{align*}
\|\bar{\yinf}\|^2_{\bar{\Hl}} &= \hat{w}^2\sum_{i=1}^\infty\|\left(S_\Sigma\hat{\uinf}\right)_i\|^2_2 = \hat{w}^2\|S_\Sigma\hat{\uinf}\|^2_2 \Rightarrow \\
\|\bar{\yinf}\|_{\bar{\Hl}}   &= \hat{w}\|S_\Sigma\hat{\uinf}\|_2 \le \hat{w}\underset{\|\hat{\uinf}\|_2\le 1}{\sup}\|S_\Sigma\hat{\uinf}\|_2\enspace.
\end{align*}
From the definition of the induced 2-norm of $\Sigma$, denoted here as $\|\Sigma\|_2$, we have that $\sup\|S_\Sigma\hat{\uinf}\|_2 = \|\Sigma\|_2\|\hat{\uinf}\|_2$ and by assuming (without loss of generality) that $\|\hat{\uinf}\|_2=\|\uinf\|_{\Hu}=1$, we end up having that
$\|\bar{\yinf}\|_{\bar{\Hl}}\le\hat{w}\|\Sigma\|_2$. Finally, the operator $\Cxinf$ is bounded by the $\mathcal{H}_\infty$ norm of the transfer matrix $H_\Sigma(z)$, or $\|\bar{\yinf}\|_{\bar{\Hl}}\le \hat{w}\mathcal{H}_\infty(\Sigma)$. Subsequently, we have from~(\ref{App.C:toprove}) that the operator $\Cinf$ is bounded by $\sigma_{\max}(C_u) + \hat{w}\mathcal{H}_\infty(\Sigma)$.
\end{IEEEproof}

\begin{myrem}
 Following the discussion from Section~\ref{sec:intro}, the weight $\hat{w}$ can be chosen to render any unstable system stable by shrinking the eigenvalues of the matrix $A$ ($\hat{A}=\hat{w}A$).  
\end{myrem}

\begin{mylem}{\label{lem:bounded_Hinv}}
 The operator $\Hinf^{-1}$ is bounded and $\| \Hinf \| \le 1/\lambda_{\min}(R)$. 
\end{mylem}
\begin{IEEEproof}
 The operator $\Hinf$ is given by $\Hinf=\Binf^{\star}\Qinf\Binf+\Rinf$; see~(\ref{operators}). Since $\Binf^{\star}\Qinf\Binf$ is positive semidefinite (i.e., $\langle \Binf^{\star}\Qinf\Binf \uinf,\uinf \rangle \ge 0$ for all $\uinf \in \mathcal{H}_{\uinf}$) and $\Rinf$ is positive definite according to Definition~\ref{def:posdef}, we conclude that $\Hinf$ is positive definite and hence has a bounded inverse.
 
 In order to compute the bound for $\Hinf^{-1}$ we observe that $\|(\Binf^{\star}\Qinf\Binf+\Rinf)^{-1}\|\le\|\Rinf^{-1}\|\le 1/\lambda_{\min}(R)$, which concludes the proof. 
\end{IEEEproof}

\section{Backtracking stepsize rule}{\label{App.D}}
Theorem~1 in~\cite{AFBS_chambolle}, proves convergence in the function values of $F = f + g$ as well as weak convergence of the iterates under the assumption of a fixed stepsize $\rho=\frac{1}{L(f)}$, with $\|\nabla f(\lambdainf_1)-\nabla f(\lambdainf_2)\|\le L(f)\|\lambdainf_1-\lambdainf_2\|$. As discussed in Section~\ref{sec:lin_sys_solve}, an upper bound for the global Lipschitz constant can be overly conservative and lead to a small stepsize. In this Appendix we briefly revise the backtracking stepsize rule that allows for local estimates of the curvature of $f$ and show that the results of Theorem~1 also hold in this case. The arguments are in line with~\cite{fista}. Since the convergence theory developed in~\cite{fista} applies to general Hilbert spaces (see Remark~2.1), the following hold in the $\lw$ space of interest.

We denote as $f(\cdot)$ the smooth and $g(\cdot)$ the nonsmooth convex functions of interest. Note that in our case $f=h^\star$ and $g=\delta_-$. 
For any $L>0$ consider the quadratic approximation of $F(\lambdainf)=f(\lambdainf)+g(\lambdainf)$ at a point $\yinf$:
\begin{equation}{\label{App::majorizer}}
Q_L(\lambdainf,\yinf) := f(\yinf) + \langle \lambdainf-\yinf, \nabla f(\yinf) \rangle + \frac{L}{2}\|\lambdainf-\yinf\|^2 + g(\lambdainf)\enspace.
\end{equation}
We also define the unique minimizer parametrized by the point $\yinf$ as
\begin{align}{\label{App::optimizer}}
p_L(\yinf) &:= \underset{\lambdainf}{\argmin}\left\{Q_L(\lambdainf,\yinf)\right\} \\
           &= \underset{\lambdainf}{\argmin}\left\{g(\lambdainf) + \frac{L}{2}\|\lambdainf-\left(\yinf-\frac{1}{L}\nabla f(\yinf)\right)\|^2\right\}\\
           &= \prox{_g}\left(\yinf-\frac{1}{L}\nabla f(\yinf)\right)\enspace,
\end{align}
which is the basic step of Algorithm~\ref{al:AFBS}, \ie, Step~3. The following holds:
\begin{mylem}[Lemma~2.3~\cite{fista}]{\label{lem:2_3_Teboulle}}
Let $\yinf\in\Hl$ and $L>0$ be such that
\[
 F(p_L(\yinf)) \le Q_L(p_L(\yinf),\yinf)\enspace.
\]
Then for any $\lambdainf\in\Hl$,
\[
 F(\lambdainf)-F(p_L(\yinf)) \ge \frac{L}{2}\|p_L(\yinf)-\yinf\|^2 + L\langle \yinf-\lambdainf,p_L(\yinf)-\yinf\rangle\enspace.
\]
\end{mylem}
The backtracking procedure as described in~\cite{fista} is as follows:
\begin{algorithm}[H]
\caption{Backtracking for stepsize computation}
\label{al:Backtracking}
\begin{algorithmic}  
\STATE 0: Take $L^0>0$, some $\eta>1$, and $\lambdainf^0\in\Hl$.
\REPEAT
\STATE 1: \scalebox{0.92}Find the smallest nonnegative integer $i^k$ such that with $\bar{L}=\eta^{i^k}L^{k-1}$
\UNTIL{$F(p_{\bar{L}}(\yinf^k))\le Q_{\bar{L}}(p_{\bar{L}}(\yinf),\yinf)$}
\end{algorithmic}
\end{algorithm}\begin{footnotesize}
\end{footnotesize}
Note that for any $L^k=\bar{L}$ generated by Algorithm~\ref{al:Backtracking} Lemma~\ref{lem:2_3_Teboulle} holds.

We also have the following instrumental Lemma from~\cite{AFBS_chambolle}:
\begin{mylem}[Lemma~1~\cite{AFBS_chambolle}]{\label{lem:1_Chambolle}}
Let $L\ge L(f)$, $\lambdainf,\yinf\in\Hl$ and $p_L(\yinf):=\prox{_g}(\yinf - \frac{1}{L}\nabla f(\yinf))$. Then for all $\lambdainf$
\[
 F(p_L(\yinf)) + \frac{L}{2}\|p_L(\yinf)-\yinf\|^2 \le F(\yinf) + \frac{L}{2}\|\lambdainf-\yinf\|^2\enspace.
\]
\end{mylem}
Lemma~\ref{lem:1_Chambolle} is crucial for proving point (i) of Theorem~\ref{theor:master}. More specifically, supposing that Lemma~\ref{lem:1_Chambolle} holds, convergence of the function values with rate $1/k^2$ can be proven under no further assumptions using Theorem~2 and Corollary~1 in~\cite{AFBS_chambolle}. We will not repeat the aforementioned theorems here due to space limitations. What we are going to show instead is that all stepsizes $\rho^k=1/L^k$ generated from Algorithm~\ref{al:Backtracking} satisfy Lemma~\ref{lem:1_Chambolle}.
\begin{mylem}{\label{lem::back_conv}}
 Consider $F=f+g$, with $f=h^\star$ and $g=\delta_-$ as defined in Section~\ref{sec:intro}. The iterates $\lambdainf^k$ generated from Algorithm~\ref{al:AFBS} with a backtracking stepsize rule generated from Algorithm~\ref{al:Backtracking} satisfy:
 \[
 F(\lambdainf^k)-F(\lambdainf^\infty) \le \frac{a^2\bar{L}}{2(k+a-1)^2}\|\lambdainf^0-\lambdainf^\infty\|^2\enspace.
 \] 
\end{mylem}
\begin{IEEEproof}
 All local Lipschitz estimates $\bar{L}$ generated from Algorithm~\ref{al:Backtracking} satisfy $F(p_{\bar{L}}(\yinf^k))\le Q_{\bar{L}}(p_{\bar{L}}(\yinf),\yinf)$ and, consequently, from Lemma~\ref{lem:2_3_Teboulle}, 
 \[
 F(\lambdainf)-F(p_{\bar{L}}(\yinf)) \ge \frac{\bar{L}}{2}\|p_{\bar{L}}(\yinf)-\yinf\|^2 + \bar{L}\langle \yinf-\lambdainf,p_{\bar{L}}(\yinf)-\yinf\rangle\enspace,
 \]
 or 
 \[
 \frac{2}{\bar{L}}\left(F(\lambdainf)-F(p_{\bar{L}}(\yinf))\right) \ge \|p_{\bar{L}}(\yinf)-\yinf\|^2 + 2\langle \yinf-\lambdainf,p_{\bar{L}}(\yinf)-\yinf\rangle\enspace.
 \]
 It follows from the Pythagorean theorem that
 \[
 \|b-a\|^2 + 2\langle b-a,a-c\rangle = \|b-c\|^2 - \|a-c\|^2\enspace,
 \]
 and hence
  \[
 \frac{2}{\bar{L}}\left(F(\lambdainf)-F(p_{\bar{L}}(\yinf))\right) \ge \|p_{\bar{L}}(\yinf)-\yinf\|^2 - \|\yinf-\lambdainf\|^2\enspace,
 \]
 from which Lemma~\ref{lem:1_Chambolle} follows. Theorem~2 and Corollary~1 of~\cite{AFBS_chambolle} then lead to the desired result.
\end{IEEEproof}

%\section{Conclusion}
%
%A conclusion section is not required. Although a conclusion may review
%the main points of the paper, do not replicate the abstract as the
%conclusion. A conclusion might elaborate on the importance of the work
%or suggest applications and extensions.

\section*{Acknowledgements}
The authors would like to thank J.H. Hours for many helpful discussions about the resulting optimization problems. \\

The research leading to these results has received funding from the European Research Council under the European Union's Seventh Framework Programme (FP/2007-2013)/ ERC Grant Agreement n. 307608. 

\ifCLASSOPTIONcaptionsoff
  \newpage
\fi  

\bibliographystyle{IEEEtran}
\bibliography{IEEEabrv,CLQR_arXiv_single}         % bib file to produce the bibliography
                                                     % with bibtex (preferred)

\end{document}